\documentclass[sn-mathphys]{sn-jnl}% Math and Physical Sciences Reference Style
\jyear{2021}%

\theoremstyle{thmstyleone}%
\theoremstyle{thmstyletwo}%

\theoremstyle{thmstylethree}%
\newtheorem{definition}{Definition}%

\usepackage{lmodern}
\usepackage[english]{babel}
\usepackage[utf8]{inputenc} %caractère accentués
\usepackage{amsmath}
\usepackage{enumitem}
\usepackage[symbol]{footmisc}
\usepackage{xparse}%
\usepackage{xkeyval}%
\usepackage{subfloat}
\usepackage{setspace} 
\usepackage{tabularx}
\usepackage{setspace} 
\usepackage{amsfonts}
\usepackage{graphicx}
\usepackage{footmisc}
\usepackage{subfig}
\usepackage{caption}
\usepackage{xcolor}
\usepackage{breqn}
\usepackage{caption}
\usepackage{float}
\usepackage{bm}
\usepackage{natbib}
\usepackage{amssymb}
\usepackage{xspace}
\usepackage{booktabs}
\usepackage{etoolbox}
\usepackage{textcomp}
\usepackage{placeins}

\newcommand{\xiz}{\widehat{x}_{i}}
\newcommand{\yiz}{\widehat{y}_{i}}
\newcommand{\xjz}{\widehat{x}_{j}}
\newcommand{\yjz}{\widehat{y}_{j}}
\newcommand{\viz}{\widehat{v}_{i}}
\newcommand{\vjz}{\widehat{v}_{j}}
\newcommand{\xijz}{\widehat{x}_{ij}}
\newcommand{\yijz}{\widehat{y}_{ij}}
\newcommand{\qi}{q_{i}}
\newcommand{\qj}{q_{j}}
\newcommand{\ti}{\theta_{i}}
\newcommand{\tj}{\theta_{j}}
\newcommand{\qmax}{\overline{q}_i}
\newcommand{\qmin}{\underline{q}_i}
\newcommand{\qub}{\overline{\bm{q}}}
\newcommand{\qlb}{\underline{\bm{q}}}

\newcommand{\pzx}{\widehat{x}_{ij}}
\newcommand{\pzy}{\widehat{y}_{ij}}
\newcommand{\tiz}{\widehat{\theta}_{i}}
\newcommand{\tjz}{\widehat{\theta}_{j}}
\newcommand{\tlb}{\underline{\bm{\theta}}}
\newcommand{\tub}{\overline{\bm{\theta}}}
\newcommand{\tmin}{\underline{\theta}_i}
\newcommand{\tmax}{\overline{\theta}_i}

\newcommand{\tm}{t^{\text{min}}_{ij}}

\newcommand{\vijx}{v^x_{ij}}
\newcommand{\vijy}{v^y_{ij}}

\renewcommand{\P}{\mathcal{P}}
\newcommand{\Po}{\mathcal{P}_0}

\newcommand{\vix}{v_{i,x}}
\newcommand{\viy}{v_{i,y}}
\newcommand{\vjx}{v_{j,x}}
\newcommand{\vjy}{v_{j,y}}

\newcommand{\Cplex}{\textsc{Cplex}\xspace}

\newcommand{\Dis}{$\mathsf{Disjunctive}$\xspace}

\newcommand{\Naive}{$\mathsf{Exact Naive}$\xspace}
\newcommand{\Greedy}{$\mathsf{Greedy Naive}$\xspace}

\newcommand{\tc}{TC}
\newcommand{\deltatc}{\Delta \tc}

\newtheorem{modl}{Model}

\newenvironment{model}{\begin{samepage}\begin{modl}}{\end{modl}\end{samepage}}

\begin{document}

\title[\small{Aircraft conflict resolution with trajectory recovery
using mixed-integer programming}]{Aircraft conflict resolution with trajectory recovery using mixed-integer programming}


\author*[1]{\fnm{Fernando H. C.} \sur{Dias}}\email{fernando.cunhadias@helsinki.fi}

\author[2]{\fnm{David} \sur{Rey}}\email{david.rey@skema.edu}

\affil*[1]{\orgdiv{Department of Computer Science}, \orgname{University of Helsinki}, \orgaddress{\postcode{00560}, \country{Finland}}\footnotemark[3]}

\affil[2]{\orgdiv{SKEMA Business School}, \orgname{Universit\'e C\^ote d'Azur}, \orgaddress{ \city{Sophia Antipolis Campus}, \country{France}}}


\abstract{To guarantee the safety of flight operations, decision-support systems for air traffic control must be able to improve the usage of airspace capacity and handle increasing demand. In this study, we address the aircraft conflict avoidance and trajectory recovery problem. The problem of finding least deviation conflict-free aircraft trajectories that guarantee the return to a target waypoint is highly complex due to the nature of the nonlinear trajectories that are sought. We present a two-stage iterative algorithm that first solves initial conflict by manipulating their speed and heading control and then identifying each aircraft's optimal time to recover its trajectory towards their nominal. The avoidance stage extends existing mixed-integer programming formulations, and for the recovery stage, we propose a novel mixed-integer formulation. We assume that speed and heading control are continuous variables for this approach while the recovery time is treated as a discrete variable. In this approach, it is shown that the trajectory recovery costs can be anticipated by inducing avoidance trajectories with higher deviation, therefore obtaining earlier recovery time within few iterations. Numerical results on benchmark conflict resolution problems show that this approach can solve instances with up to 30 aircraft within 10 minutes.}


\footnotetext[3]{This work was partially funded by the European Research Council (ERC) under the European Union's Horizon 2020 research and innovation programme (grant agreement No.~851093, SAFEBIO)}


%\footnotetext[3]{This work was partially funded by the European Research Council (ERC) under the European Union's Horizon 2020 research and innovation programme (grant agreement No.~851093, SAFEBIO).\label{mark}}
\maketitle

\section{Introduction}\label{sec1}

Airspace usage has seen an increase in demand throughout the last decades (except during specific periods such as post 09/11 attacks and the COVID-19 pandemic). Aligned with this trend, the limited airspace has put the current air traffic management (ATM) system under intense pressure. This is reflected in the increasing amount of control necessary to guarantee safety, which is a crucial aspect of air traffic control (ATC). With the advance of unmanned aerial systems and urban air mobility, denser and congested traffic configurations are also expected. This configuration may lead to impairment of aircraft safety.
Nevertheless, state-of-the-art methods for aircraft traffic control are reaching their limits, and new approaches, including more automation, have recently received significant attention in the field \citep{durand2009ant, vela2010near}. Introducing automation within ATC systems has the potential to reduce controller workload and improve airspace capacity \citep{rey2016subliminal}. Conflict detection and resolution (CDR) is vital for the workload model of air traffic controllers. Hence, this calls for advanced conflict algorithms capable of acting as efficient decision-support tools. 

The aircraft conflict avoidance and resolution problem can be formulated as an optimization problem in which the goal is to find conflict-free trajectories while minimizing a cost function, e.g. the deviation to the original flight path. Many strategies have been proposed to address this problem based on the type of manoeuvres issued to aircraft: speed, heading or altitude control. These strategies can be applied separately or in combination. In this study, we refer to conflict avoidance all conflict resolution approaches that do not consider the aircraft trajectory recovery problem. Conflict resolution using global optimization has recently received growing attention due to its ability to provide optimal solutions that consider all traffic within an airspace region. Following, we review the literature on conflict avoidance in Section \ref{litconf} before focusing on efforts to develop approaches for trajectory recovery in Section \ref{litrec} We outline our contributions relative to the field in Section \ref{litconf}.

\subsection{Conflict Avoidance}
\label{litconf}

One of the first global optimization approaches for aircraft conflict resolution was introduced by \cite{pallottino2002conflict} which proposed two formulations: one focusing on speed control and another focusing on heading control, and both minimize overall flight time. In the proposed MIP (Mixed Integer Problem) formulation for conflict resolution with speed control, the authors derived linear pairwise aircraft separation constraints based on the geometric construction introduced by \cite{bilimoria2000geometric}. These separation conditions are obtained by projecting the shadow of an intruder aircraft onto the trajectory of a reference aircraft. \cite{frazzoli2001resolution} was the first to observe that this geometric construction provided a basis to characterize the set of aircraft pairwise conflict-free trajectories via linear half-planes in the relative velocity (speed and heading) plane. The authors introduced a nonconvex formulation for the conflict resolution problem with speed and heading control and proposed a convex relaxation based on semi-definite programming and a heuristic algorithm to find feasible solutions on problems with up to 10 aircraft.

Subsequent approaches proposed speed control and altitude level-assignment to minimize fuel consumption by metering aircraft at conflict points \citep{vela2010near}. In \cite{vela2009two}, the authors proposed a two-stage stochastic optimization model accounting for wind uncertainty and using speed control. Multi-objective optimization formulations attempting to balance flight deviation with the total number of manoeuvres (velocity, heading or altitude change), building on the work of \cite{pallottino2002conflict} were proposed by \citep{alonso2011collision,alonso2014exact}. Subliminal speed control methods which focus on speed control only for conflict resolution has also proven to be efficient and with low impact in terms of deviation and fuel consumption, although they may fail to resolve all conflicts \citep{rey2016subliminal,cafieri2017maximizing}. 

More recently, nonlinear global optimization approaches received increasing attention in the literature. \cite{omer2013hybridization} proposed a hybrid algorithm that uses the optimal solution of a MILP (Mixed Integer Linear Problem) as the starting point for solving a nonlinear formulation of the same problem. \cite{cafieri2017maximizing} proposed an MINLP (Mixed Integer Non-Linear Problem) approach for conflict resolution with speed control only which highlights that subliminal speed control alone may not be sufficient to resolve all conflicts in dense traffic scenarios. Using a similar framework, \cite{cafieri2017mixed} presented a two-step approach where a maximum number of conflicts are first solved using speed control only, and outstanding conflicts are solved by heading control. \cite{cerulli2020detecting} proposed a formulation based on bi-level optimization with multiple follower problems, each of which represents a two-aircraft separation problem. The authors presented two formulations, one using speed control and another using heading control. A cut generation algorithm is proposed to solve the corresponding bi-level optimization problems. More recently, \cite{rey2017complex} proposed a complex number formulation for speed and heading control without any form of discretization that was further expanded by \cite{dias2021disjunctive} in which an exact constraint generation algorithm was proposed. A systematic review of mathematical programming methods in ATC can be found in \cite{pelegrin2020airspace}.

\subsection{Trajectory Recovery}
\label{litrec}

In mathematical programming approaches, aircraft trajectory is usually divided into two stages. Those stages usually are called avoidance (or action) which encapsulates the manoeuvres necessary to avoid any conflict, while the trajectory recovery corresponds to the manoeuvres necessary to restore the aircraft to its original trajectory. Different approaches (such as genetic programming and heuristics) might be considered trajectory recovery as part of conflict resolution. However, this is usually separated in mathematical programming due to its non-trivial aspects, bringing many non-convexity and non-linearity.

Despite their potential effectiveness, most efforts in conflict resolution have focused on ensuring conflict avoidance but have overlooked the costs and mechanisms for modelling aircraft's recovery to their original trajectory. This may be critical when conflict resolution is performed using heading control which may significantly deviate the aircraft from their initial trajectory, thus possibly increasing flight operating costs. Aircraft conflict resolution with trajectory recovery is the problem of finding conflict-free trajectories which ensure that aircraft recover their initial flight path upon completion of the manoeuvres performed. This problem has received very little attention in the literature due to its challenging nature. Meta-heuristics such as genetic algorithms exploited the fact that formulating an analytical trajectory is quite challenging (due to trigonometric functions) and expressed the aircraft path by discretizing it. In \cite{durand1997optimal}, the aircraft path is divided into different segments based on possible conflict points and turning points, thus predefining a set of eligible trajectories. An ant colony algorithm was also proposed by \cite{durand2009ant} where the authors determined the target point for each aircraft trajectory and discretized its trajectory using a specific timestep and in \cite{durand1996automatic}, where another formulation is attempted with discretized space and time. In order to generate a smooth aircraft path and derive conflict-free solutions, genetic algorithms were implemented. \cite{dougui2013light} proposed a model which uses an analogy with light propagation theory to create conflict-free aircraft trajectories with recovery. \cite{peyronne2015solving} proposed a B-splines model which uses way-points of a given trajectory to design conflict-free trajectories with recovery. In \cite{omer2015space}, the author proposed a formulation providing parallel trajectory recovery while minimizing fuel consumption and delays. In this model, aircraft are assumed to perform a preventive manoeuvre before their intersection with other trajectories, and the formulation aims to find trajectories that are parallel to the aircraft initial's trajectory. Heading angles are discretized, and the optimization controls both aircraft heading and recovery time. Recently, \cite{lehouillier2017two} proposed a manoeuvre-discretized model in which predefined sets of manoeuvres are available for aircraft, and a clique-based formulation is proposed to find the optimal combination of conflict-free manoeuvres. 

This literature review highlights that despite recent improvements in computational optimization, there remain significant open challenges in designing scalable global optimization approaches for conflict resolution in air traffic control, especially on incorporating recovery in a scalable and effective way. Only a few methods for aircraft trajectory recovery have been proposed, and they typically assume a simplified trajectory design (as expressed beforehand). This can be explained by the complexity and nonlinear aspects present in the trajectory recovery formulation, especially using mathematical programming. Approaches that have jointly addressed the avoidance and recovery problems often assume the existence of complete manoeuvre sets before optimization or any other forms of variable discretization or simplification.

\subsection{Our Contributions}
\label{contr}

As highlighted in our literature review, most of the state-of-art models using mathematical programming focus on conflict avoidance only. Specifically, they combine speed and heading control manoeuvres in deterministic approaches. In this study, we present a new two-stage algorithm for aircraft conflict resolution with trajectory recovery. In this approach, the speed and heading of aircraft are first optimized to avoid conflicts while minimizing the deviation from their initial trajectories. Then, in a second stage, aircraft trajectories are modified to recover a target position on the aircraft's initial trajectories. These two stages are incorporated into an iterative algorithm, where the recovery cost is projected into the avoidance stage. Hence, this algorithm aims to obtain a non-trivial solution where the overall cost of aircraft manoeuvres is minimized. The main contributions of this study relative to the literature are i) the extension of state-of-the-art mixed-integer formulations for conflict avoidance to include on/off variables and constraints to select aircraft that are being controlled in the conflict resolution process; ii) the design of a penalty-based algorithm that iterates the conflict-avoidance and trajectory recovery stages to find non-trivial solutions to the problem at hand, and iii) extensive numerical experiments on benchmarking problems from the literature that demonstrate the benefits of the proposed approach. 

The paper is organized as follows: in Section \ref{two_stage}, the two-stage algorithm is detailed, starting with the conflict avoidance formulation and then the trajectory recovery portion. In Section \ref{num_contin}, the experimental framework, the results and discussions are presented. In Section \ref{con-rec}, the conclusion and future research directions are provided. 

\section{Conflict Resolution with Trajectory Recovery}
\label{two_stage}

In this section, we present a novel approach for conflict resolution with trajectory recovery that is based on decomposing the problem into two stages. The proposed approach aims to repeatedly solve these two sub-problems, i.e. conflict avoidance and trajectory recovery using a penalty-based mechanism to anticipate recovery costs in the avoidance stage.

\subsection{Preliminaries}

In aircraft trajectory formulation, we assume that aircraft current and target positions are known, and that aircraft are initially not in conflict. This sets the context of the optimization problem of interest: given a set of aircraft with known current and target positions, find least-deviating conflict-free trajectories for all aircraft, such that aircraft may safely reach their target destination. To address this problem, we propose decomposing the trajectory optimization problem into two stages: 1) conflict avoidance and 2) trajectory recovery. The first stage focuses on controlling aircraft heading and speed to avoid all conflicts, while the second stage focuses on calculating the optimal time for aircraft to start safely recovering towards their target position. We focus on the two-dimensional conflict resolution problem and only consider horizontal aircraft manoeuvres for brevity. The extension to the vertical case can be addressed by incorporating flight level change manoeuvres in the conflict avoidance stage \citep{dias2021disjunctive} and ensuring safe recovery to aircraft target flight level and position. We leave this extension for future research. 

\begin{figure}[htp]
    \centering
    \includegraphics[width=1.0\textwidth]{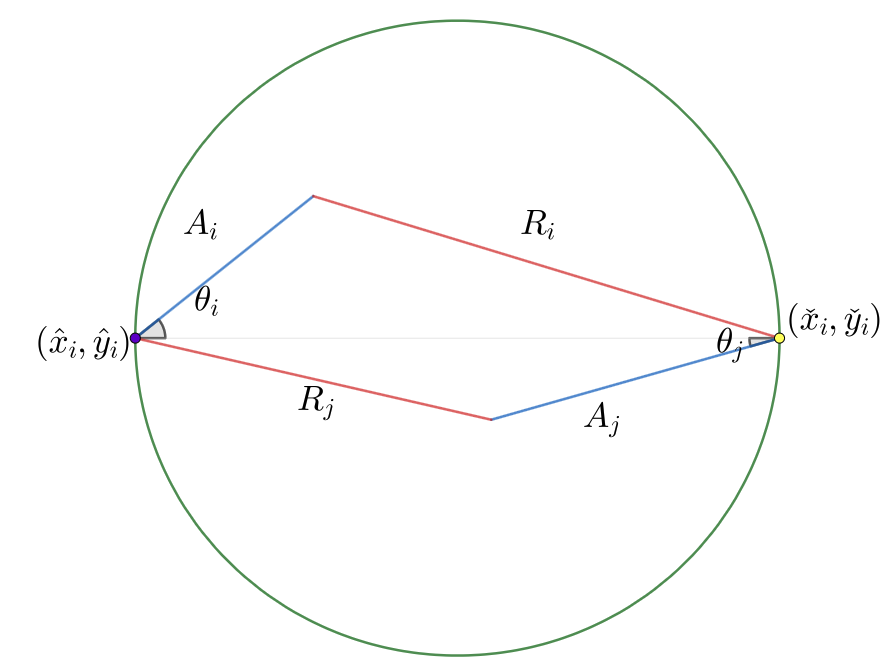}
    \caption{Example of a two-aircraft conflict resolution problem. If both aircraft followed their nominal trajectory (in grey), aircraft $i$ is in conflict with aircraft $j$. In order to avoid conflict, conflict avoidance manoeuvres have to be taken. For aircraft $i$, it deviates from its nominal manoeuvre by following the avoidance trajectory ($A_i$). After a certain amount of time, the aircraft can safely return towards its target point ($\check{x}_i$,$\check{y}_i$) via the trajectory recovery ($R_i$).}
    \label{fig:traj}
\end{figure}

The avoidance stage aims to determine the optimal variation in speed and angle for each aircraft to avoid conflicts. In the recovery stage, it is necessary to calculate and identify the manoeuvres required to recover aircraft initial trajectories. In order to impose airspace safety, it is necessary to impose the separation conditions in the recovery stage as well. In mathematical programming, recovery approaches are challenging for several reasons. Assuming that aircraft trajectories can be modelled as linear and that uniform motion laws hold, the conflict separation conditions are created via Euclidean distance. This leads to equations using quadratic and trigonometric elements. This creates difficulty in organizing and developing an optimization model using those equations. Thus, such formulations may not scale easily and be used in larger scenarios. Second, Euclidean distance creates quadratic components that are also nonlinear and nonconvex concerning decision variables and cannot be easily solved. As presented by \cite{pallottino2002conflict}, \cite{alonso2011collision} and \cite{rey2017complex} nonlinear constraints can be further simplified into a set of integer-linear components with regards to aircraft velocity. However, these formulations only focus on separating aircraft to solve existing conflicts and ignore the process of planning aircraft recovery to a target destination. As described Fig \ref{fig:traj}, trajectories such as the one highlighted in the blue and red segments are an example of conflict avoidance and trajectory recovery, respectively. At first, the aircraft manoeuvre as the blue segment represents. Another manoeuvre is performed to redirect the aircraft towards its target point (in yellow). In this study, we focus on developing a global optimization approach to construct piecewise linear conflict-free trajectories (avoidance and recovery) for a set of aircraft to minimize the total deviation from their initial original trajectories.

\subsection{Conflict Avoidance}
Consider a set of aircraft $\mathcal{A}$ sharing the same flight level. For each aircraft $i \in \mathcal{A}$, assuming uniform motion laws apply, its position is: $p_i(t) = [x_i(t) = \widehat{x}_i + q_i\hat{v}_i\cos (\widehat{\theta}_{i} + \theta_i)t, y_i(t) = \widehat{y}_i + q_i\hat{v}_i\sin(\widehat{\theta}_{i} + \theta_i)t]^\top$ in which $v_i$ is the speed, $\widehat{x}_i$ and $\widehat{y}_i$ are the initial coordinates of $i$ at the time of optimization, $\widehat{\theta}_{i}$ is its initial heading angle, $\theta_i$ is its deviation angle and $q_i$ is the speed deviation. The relative velocity vector of $i$ and $j$, denoted $v_{ij}$, can be expressed as:

\begin{equation}
    v_{ij} = [\vijx,\vijy]^\top = [\vix - \vjx,\viy - \vjy]^\top,
\end{equation}

where:
\begin{align}
    \vix = q_i\hat{v}_i\sin(\widehat{\theta}_{i} + \theta_i), \\
    \viy = q_i\hat{v}_i\cos(\widehat{\theta}_{i} + \theta_i). 
\end{align}

Incorporating these elements into the equation of motion gives:
\begin{subequations}
\begin{align}
& {x}_{i}(t) = \xiz + \vix t,\\
& {y}_{i}(t) = \yiz + \viy t.
\end{align}\label{randommotion}
\end{subequations}

The relative position of aircraft $i$ and $j$ at time $t$ can be represented as $p_{ij}(t) = p_{i}(t) - p_{j}(t)$. Let $d$ be the horizontal separation norm, typically $d = 5$ NM. Two aircraft $i,j \in \mathcal{A}$ are horizontally separated if and only if: $||p_{ij}(t)|| \ge d$, for all $t \ge 0$.

Let $\P$ be the set by the pairs of aircraft, i.e. $\P = \{i \in \mathcal{A}, j \in \mathcal{A}, i <j\}$. For each pair $(i,j) \in \P$, the relative position vector ${p}_{ij}(t)$ is:
\begin{equation}
p_{ij}(t) = [x_{ij}(t),y_{ij}(t)]^\top,
\end{equation}

and the relative velocity vector ${v}_{ij} = [\vijx,\vijy]^\top,$ is:
\begin{subequations}\label{eq:v}
\begin{align}
& \vijx = \qi \viz\cos(\tiz + \ti) - \qj \vjz\cos(\tjz + \tj), \\
& \vijy = \qi \viz\sin(\tiz + \ti) - \qj \vjz\sin(\tjz + \tj).
\end{align}
\end{subequations}

Imposing the separation condition, gives for each pair $(i,j) \in \P$:
\begin{equation}\label{eq:sepcond}
||{p}_{ij}(t)|| \geq d \Leftrightarrow \sqrt{({x}_i(t) - y_j(t))^2 + ({y}_i(t)- {y}_j(t))^2} \geq d, \quad \forall t \ge 0.
\end{equation}

Let $\xijz = \xiz - \xjz$ and $\yijz = \yiz - \yjz$. Squaring both sides in Eq. \eqref{eq:sepcond}, we obtain:
\begin{equation}
{f}_{ij}(t) \equiv ((\vijx)^2 + (\vijy)^2)t^2 + (2\vijx\xijz + 2\vijy\yijz)t + \xijz^2 + \yijz^2- d^2 \geq 0.
\end{equation}

The function ${f}_{ij}(t)$ is a second-order polynomial in $t$ which is minimal for ${f}^\prime_{ij}(t) = 0$:
\begin{equation}\label{trobust}
{f}^\prime_{ij}(t) = 0 \Rightarrow \tm \equiv -\frac{\xijz\vijx + \yijz\vijy}{(\vijx)^2 + (\vijy)^2}.
\end{equation}

The time instant $\tm$ represents the time of minimal separation of aircraft $i$ and $j$. As noted in several studies \citep{cafieri2017maximizing,cafieri2017mixed,rey2017complex}, if $\tm \leq 0$ then aircraft $i$ and $j$ are diverging and, assuming aircraft are separated at $t=0$, they are thus separating for any $t \geq 0$. Further, substituting $\tm$ in ${f}_{ij}(t)$ yields:
\begin{equation}\label{grobust}
{g}_{ij}(\vijx,\vijy) \equiv {f}_{ij}(\tm) = (\vijy)^2(\xijz^2 - d^2) + (\vijx)^2(\yijz^2 - d^2) - (2 \xijz\yijz\vijx\vijy).
\end{equation}

Hence, if $g_{ij}(\vijx,\vijy) \geq 0$, then aircraft $i$ and $j$ are separated. Writing  the terms $g_{ij}$ and $\tm$ as functions of the relative velocity variables $\vijx$ and $\vijy$, we obtain the following disjunctive pairwise aircraft separation conditions:
\begin{equation}\label{eq:sepconditions}
||{p}_{ij}(t)|| \geq d, \forall t \geq 0 \Leftrightarrow {g}_{ij}(\vijx,\vijy) \geq 0 \vee \tm(\vijx,\vijy) \leq 0.
\end{equation}

The separation condition \eqref{eq:sepconditions} can be further linearized following the approach described by \cite{rey2017complex} and \cite{dias2021disjunctive}. By alternatively fixing variables $\vijx$ and $\vijy$ and solving the resulting quadratic equations, we can obtain the solution for $g_{ij}(\vijx,\vijy) = 0$. By isolating each variable, we obtain the discriminants:
\begin{equation}\label{discriminant}
\begin{cases}
\Delta_{\vijx} = 4d^2(\vijy)^2(\xijz^2 + \yijz^2 - d^2),\\
\Delta_{\vijy} = 4d^2(\vijx)^2(\xijz^2 + \yijz^2 - d^2).
\end{cases}
\end{equation}

Assuming aircraft are initially separated, then $\xijz^2 + \yijz^2 - d^2 \geq 0$ holds and thus the discriminants are positive, and the roots of equation $g(\vijx,\vijy) = 0$ are the lines defined by the system of equations:
\begin{subequations}
\begin{align}
(\hat{y}_{ij}^2 -d^2)\vijx - (\hat{x}_{ij}\hat{y}_{ij} + d\sqrt{\hat{x}_{ij}^2 + \hat{y}_{ij}^2 - d^2})\vijy = 0, \\
(\hat{y}_{ij}^2 -d^2)\vijx - (\hat{x}_{ij}\hat{y}_{ij} - d\sqrt{\hat{x}_{ij}^2 + \hat{y}_{ij}^2 - d^2})\vijy = 0, \\ 
(\hat{x}_{ij}^2 -d^2)\vijy - (\hat{x}_{ij}\hat{y}_{ij} + d\sqrt{\hat{x}_{ij}^2 + \hat{y}_{ij}^2 - d^2})\vijx = 0, \\
(\hat{x}_{ij}^2 -d^2)\vijy- (\hat{x}_{ij}\hat{y}_{ij} - d\sqrt{\hat{x}_{ij}^2 + \hat{y}_{ij}^2 - d^2})\vijx = 0. 
\end{align}
\label{robustmotion}
\end{subequations}

\begin{figure}[!h]
\centering
\subfloat[The hashed orange region represents $g(\vijx,\vijy) \geq 0$. The hashed blue half-plane represents diverging trajectories, i.e. $\tm(\vijx,\vijy) \leq 0$.]{\includegraphics[width=0.75\linewidth]{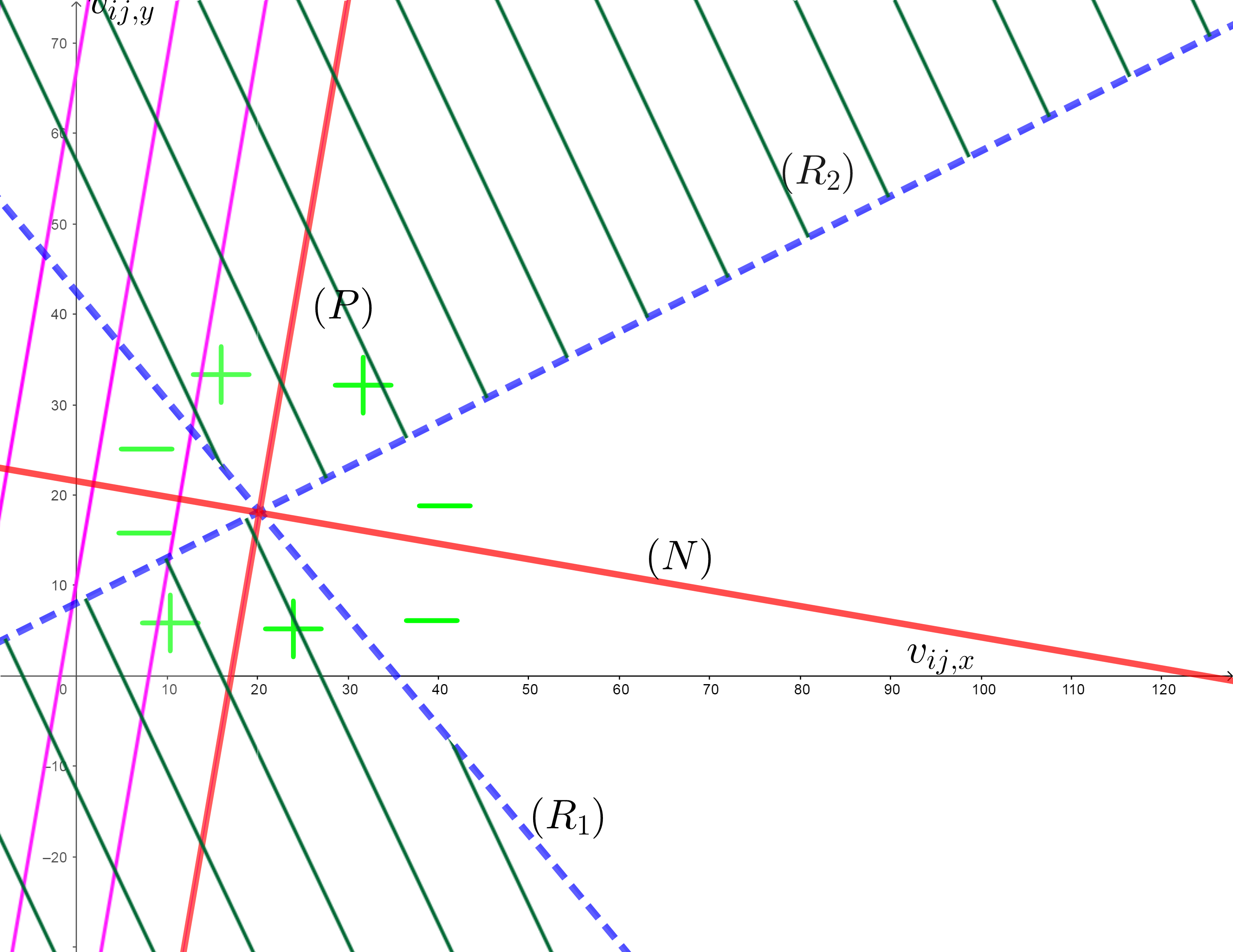}\label{fig:gplotrobust}}
\vfill
\subfloat[The inner box with black lines corresponds to the velocity bounds $\mathcal{B}$ in the deterministic scenario. The region is hashed in red corresponds to the conflict region $\mathcal{C}$. If relative velocity $\mathcal{B}$ intersects with the conflict region $\mathcal{B}$, then there exists a risk of conflict.] {\includegraphics[width=0.75\textwidth]{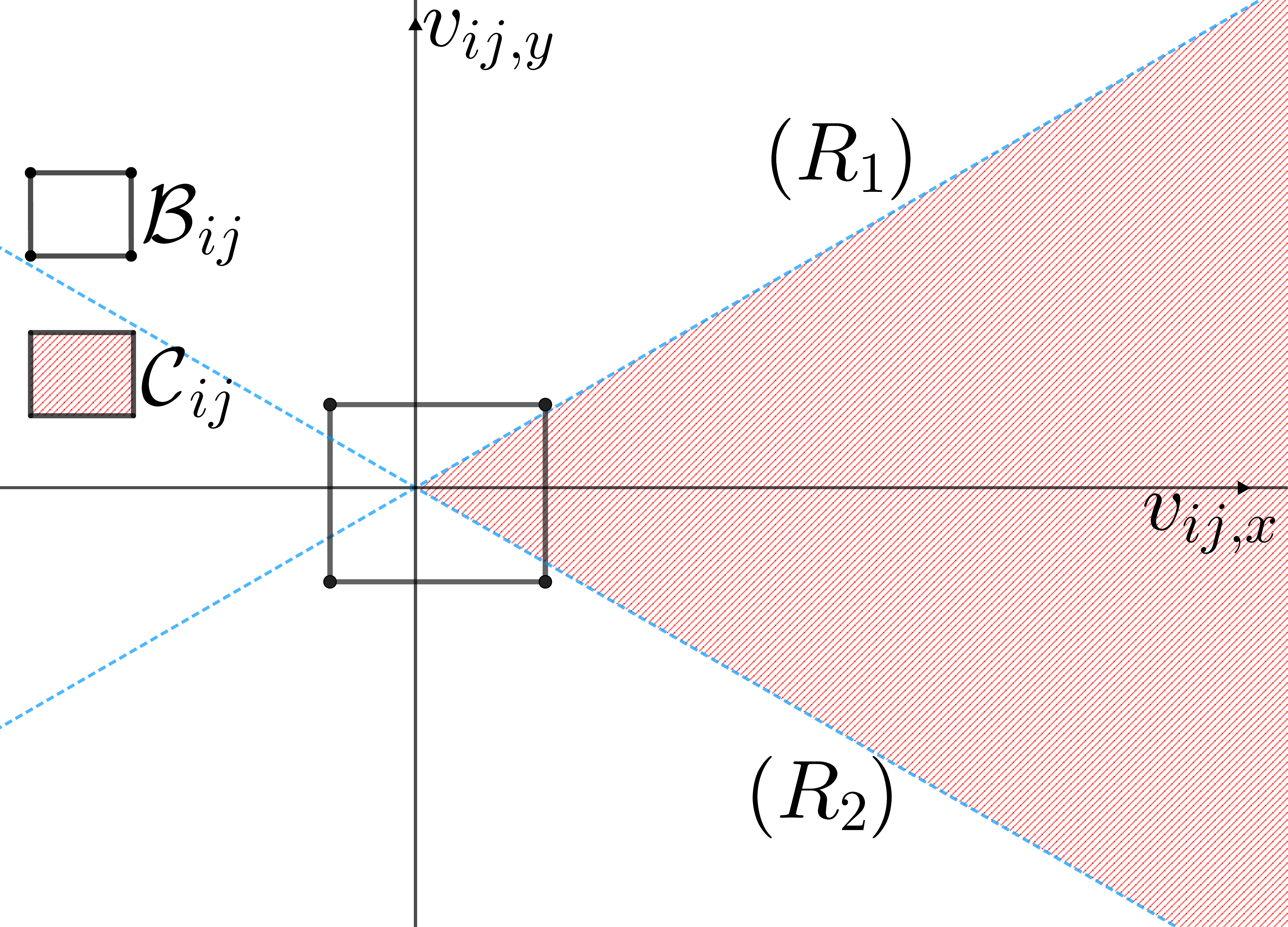}\label{conflictRegion}}
\caption{Illustration of a two-aircraft conflict in the plane $\{(\vijx,\vijy) \in \mathbb{R}^2\}$. The red lines represent the lines $P$ and $N$. The dashed blue lines correspond to the linear equations $R_1$ and $R_2$ that are the roots of $g(\vijx,\vijy) = 0$. The sign of $g(\vijx,\vijy)$ is shown by the + and - pink symbols.}
\end{figure}

If all coefficients in Eq. \eqref{robustmotion} are non-zero, they define two lines, denoted $R_1$ and $R_2$, in the plane $\{(\vijx,\vijy) \in \mathbb{R}^2\}$ and the sign of $g_{ij}(\vijx,\vijy)$ can be characterized based on the position of $(\vijx,\vijy)$ relative to these lines (see Figure \eqref{fig:gplotrobust}). Recall that according to Eq. \eqref{trobust}, the sign of the dot product~$\hat{p}_{ij} \cdot v_{ij}$ indicates aircraft convergence or divergence. Let \eqref{plane} be the equation of the line corresponding to the dot product~$\hat{p}_{ij} \cdot v_{ij}$. 
\begin{equation}\label{plane}\tag{$P$}
\vijx\xijz + \vijy\yijz = 0.
\end{equation}

The line defined by \eqref{plane} splits the plane $\{(\vijx,\vijy) \in \mathbb{R}^2\}$ in two half-planes, each of which representing converging and diverging trajectories, respectively. Consider the line normal to \eqref{plane}, denoted \eqref{normalplane}:
\begin{equation}\label{normalplane}\tag{$N$}
\vijy\hat{x}_{ij} - \vijx\hat{y}_{ij} = 0.
\end{equation}

Recall that any point $(\vijx,\vijy)$ such that $\tm \leq 0$ or $g_{ij}(\vijx,\vijy) \geq 0$ corresponds to a pair of conflict-free trajectories. Hence, the conflict-free region is nonconvex. According to the speed profile, we can divide the solution space into two regions: the region hashed in red corresponds to the conflict region $\mathcal{C}$ while the complement of this region represents conflict-free trajectories. We recall the formal definition of the conflict region. 

\begin{definition}[Conflict region, \cite{dias2021disjunctive}]
\label{cr}
Consider a pair of aircraft $(i,j) \in \P$. Let $\mathcal{C}$ be the subset of $\mathbb{R}^2$ defined as:
\begin{equation}
\mathcal{C} = \left\{(\vijx,\vijy) \in \mathbb{R}^2 : \vijx \gamma_{ij}^l - \vijy \phi_{ij}^l \geq 0 \land \vijx \gamma_{ij}^u - \vijy \phi_{ij}^u \leq 0 \right\}.
\end{equation}
$\mathcal{C}$ is the conflict region of $(i,j) \in \P$.
\end{definition}

The conflict region of pair of aircraft represents the set of relative velocity vectors $(\vijx,\vijy)$ which corresponds to conflicts. As depicted in Figure \ref{conflictRegion}, the feasible region is nonconvex, which can be solved via a disjunctive formulation. However, as presented in \cite{rey2017complex} and \cite{dias2021disjunctive}, the set of linear equations described by Eqs. \ref{eq:zd} is equivalent to \eqref{grobust} as detailed by Theorem 1 in \cite{dias2021disjunctive}. In each convex sub-region, the lines delineate the conflict-free region. The expressions of these lines depend on initial aircraft positions, i.e. $\xijz$, $\yijz$. Integer-linear separation conditions with regards to aircraft velocity components can be derived as follows, and we model this disjunction using the variable $z_{ij} \in \{0,1\}$ defined as:
\begin{subequations}
\label{eq:zd}
\begin{align}
\vijy \xijz - \vijx \yijz \leq 0, &\quad \text{ if } z_{ij} = 1, \quad \forall (i,j) \in \P,\\
\vijy \xijz - \vijy \yijz \geq 0, &\quad \text{ if } z_{ij} = 0, \quad \forall (i,j) \in \P,\\
\vijy \gamma_{ij}^l - \vijx \phi_{ij}^l \leq 0, &\quad \text{ if } z_{ij} = 1, \quad \forall (i,j) \in \P,\\
\vijy \gamma_{ij}^u - \vijx \phi_{ij}^u \geq 0, &\quad \text{ if } z_{ij} = 0, \quad \forall (i,j) \in \P.
\end{align}
\end{subequations}

In conflict avoidance problems, a common objective is to minimize the combined deviations of all aircraft. This may lead several aircraft to perform minimal conflict avoidance manoeuvres, which may not be desirable from an operational perspective. As observed in \cite{dias2020two}, even small deviations may result in costly recovery manoeuvres. In order to generate trajectories where fewer aircraft are controlled and those that are controlled have a reduced total cost, an additional binary variable is introduced. This variable determines whether an aircraft is controlled or not. In this formulation, the avoidance stage has as the objective to minimize the total deviation and the number of aircraft that are controlled. Let $f_i \in \{0,1\}$ for $i \in \mathcal{A}$ represent the control variable: $f_i = 1$ represents the case where aircraft $i$ modifies its speed and/or heading; and $f_i = 0 $ represents the case where aircraft $i$ does not perform any conflict avoidance manoeuvre.

Based on the binary variable $f_i$, if aircraft $i$ does not perform any manoeuvre, its speed and heading control should remain unchanged. This can be expressed as:

\begin{subequations}\label{eq:fqtheta}
\begin{align}
&\qmin f_i + (1 - f_i) \leq q_i \leq \qmax f_i + (1 - f_i) , && \forall i \in \mathcal{A}, \label{eq:fq}\\
&\underline{\theta}_i f_i \leq \theta_i \leq \overline{\theta}_i f_i, && \forall i \in \mathcal{A}.  \label{eq:ftheta}
\end{align}
\end{subequations}

Considering this control variable, if a pair of aircraft is in conflict, at least one of these aircraft must perform manoeuvres. Let $\Po$ be the set the aircraft that are initially in conflict such as that $\Po = \{(i,j) \in \P : {g}_{ij}^0((\vijx)^0,(\vijy)^0) \geq 0 \vee (\tm)^0((\vijx)^0,(\vijy)^0) \leq 0 \}$. Recall that a pair $(i,j)$ of aircraft is initially in conflict if and only if condition \eqref{eq:sepconditions} is not satisfied. Therefore, for all pairs of aircraft that are initially in the conflict, the following cuts are valid inequalities: 

\begin{equation}
    \label{fcons}
    f_i + f_j \geq 1, \quad \forall (i,j) \in P_0
\end{equation}

These cuts are added to strengthen the conflict resolution formulations by observing that at least one of them must perform an avoidance manoeuvre for any pair of aircraft initially in conflict. 

To incorporate the impact of controlling aircraft in the objective function, we propose to minimize a weighted sum of two terms: a fixed cost linked to control variables ($f_i$) and a variable cost linked to trajectory deviation variables ($q_i$, $\theta_i$). We denote $\lambda_f$ the weight representing the fixed cost of controlling an aircraft, and we denote $w$ in $[0,1]$ the weight used to capture the trade-off between heading and speed control deviations.

\begin{subequations}
\begin{align}
& \min \sum_{i \in \mathcal{A}} w\theta_i^2 + (1 - w)(1 - q_i)^2 + \lambda_f f_i, \nonumber \\
\end{align}
\label{obj}
\end{subequations}

For each aircraft $i \in \mathcal{A}$, we assume that the speed rate variable is lower bounded by $\qmin$ and upper bounded by $\qmax$, i.e.:
\begin{equation}\label{eq:qbound}
\qmin \leq \qi \leq \qmax, \qquad \forall i \in \mathcal{A}.
\end{equation}

We assume that the heading deviation is lower bounded by $\tmin$ and upper bounded by $\tmax$, i.e.:
\begin{equation}\label{eq:tbound}
\tmin \leq \ti \leq \tmax,  \qquad \forall i \in \mathcal{A}.
\end{equation}

To derive lower and upper bounds on relative velocity components $\vijx$ and $\vijy$, we re-arrange Eq. \eqref{eq:v} using trigonometric identities:
\begin{subequations}\label{eq:v2}
\begin{align}
& \vijx = \qi \viz\cos(\tiz)\cos(\ti) - \qi \viz\sin(\tiz)\sin(\ti) - \nonumber\\
& \qj \vjz\cos(\tjz)\cos(\tj) + \qj \vjz\sin(\tjz)\sin(\tj), \\
& \vijy = \qi \viz\sin(\tiz)\cos(\ti) + \qi \viz\cos(\tiz)\sin(\ti) - \nonumber \\
& \qj \vjz\sin(\tjz)\cos(\tj) - \qj \vjz\cos(\tjz)\sin(\tj).
\end{align}
\end{subequations}

Let $\underline{v}_{ij,x},\overline{v}_{ij,x}$ and $\underline{v}_{ij,y},\overline{v}_{ij,y}$ be the lower and upper bounds for $\vijx$ and $\vijy$, respectively. These bounds can be determined using Eq. \eqref{eq:v2} and the bounds on speed and heading control provided in Eqs. \eqref{eq:qbound} and \eqref{eq:tbound}. The derived bounds on the relative velocity components can be used to define a box in the plane $\{(\vijx,\vijy) \in \mathbb{R}^2\}$.

\begin{definition}[Relative velocity box]
\label{box}
Consider a pair of aircraft $(i,j) \in \P$. Let $\mathcal{B}$ be the  subset of $\mathbb{R}^2$ defined as
\begin{equation}
\B \equiv \left\{(\vijx,\vijy) \in \mathbb{R}^2 : \underline{v}_{ij,x}\leq \vijx \leq \overline{v}_{ij,x}, \underline{v}_{ij,y} \leq \vijy \leq \overline{v}_{ij,y}\right\}.
\end{equation}
$\mathcal{B}$ is the relative velocity box of $(i,j) \in \P$.
\end{definition}

The relative velocity box $\mathcal{B}$ characterizes all possible trajectories for the pair $(i,j) \in \P$ based on the available 2D deconfliction resources, i.e. speed and heading controls. To characterize the set of conflict-free trajectories of a pair of aircraft $(i,j) \in \P$, we compare the relative position of the relative velocity box $\mathcal{B}$ with the conflict region of this pair of aircraft. 

The conflict avoidance stage can be formulated exactly as presented in \cite{dias2021disjunctive} incorporating the penalty control variable $f_i$ as described in Eq. \eqref{eq:fqtheta} and constraint \eqref{fcons}.

\begin{model}[MINLP Formulation For Conflict Avoidance]
\label{mod_disj_rec}
\allowdisplaybreaks
\begin{subequations}
\begin{align}
&\emph{Minimise } \quad \sum_{i \in \mathcal{A}} w\theta_i^2 + (1 - w)(1 - q_i)^2 + \lambda_f f_i, \nonumber \\
&\emph{Subject to:}  && \nonumber \\
& \vijx = \qi \viz\cos(\tiz + \ti) - \qj \vjz\cos(\tjz + \tj), && \forall (i,j) \in \P,\\
& \vijy = \qi \viz\sin(\tiz + \ti) - \qj \vjz\sin(\tjz + \tj), && \forall (i,j) \in \P,\\
& \vijy\hat{x}_{ij} - \vijx\hat{y}_{ij} \leq 0, \quad \text{ if } z_{ij} = 1, && \forall (i,j) \in \P,\\
& \vijy\hat{x}_{ij} - \vijx\hat{y}_{ij} \geq 0, \quad \text{ if } z_{ij} = 0, && \forall (i,j) \in \P,\\
& \vijy \gamma_{ij}^l - \vijx \phi_{ij}^l \leq 0, \quad \text{ if } z_{ij} =1, && \forall (i,j) \in \P, \\
& \vijy \gamma_{ij}^u - \vijx \phi_{ij}^u \geq 0, \quad \text{ if } z_{ij} =0, && \forall (i,j) \in \P, \\
& f_i + f_j \geq 1, &&\forall (i,j) \in \Po, \\ 
& \qmin \leq \qi \leq \qmax, && \forall i \in \mathcal{A}, \\
& \tmin \leq \ti \leq \tmax, && \forall i \in \mathcal{A}, \\
& \vijx,\vijy \in \mathcal{B}, &&\forall (i,j) \in \P, \\
& f_i \in  \{0,1\}, &&\forall i \in \mathcal{A},\\
& z_{ij} \in \{0,1\}, &&\forall (i,j) \in \P.
\end{align}
\end{subequations}
\end{model}

This formulation is nonconvex due to the velocity constraint, which is nonconvex quadratic (Eq. \eqref{eq:v}). This results in a MINLP formulation which is challenging to solve and does not scale easily. Coefficients $\gamma_{ij}^l$, $\phi_{ij}^l$ and $\gamma_{ij}^u$, $\phi_{ij}^u$ (present in \eqref{eq:zd}) can be pre-processed based on the sign of $\pzx$ and $\pzy$. Finally, $\mathcal{B}$ represents the bounds for the velocity variables $(\vijx,\vijy)$.

\subsection{Trajectory Recovery}
\label{twostage}

\begin{figure}[htp]
    \centering
    \includegraphics[width=1.0\linewidth]{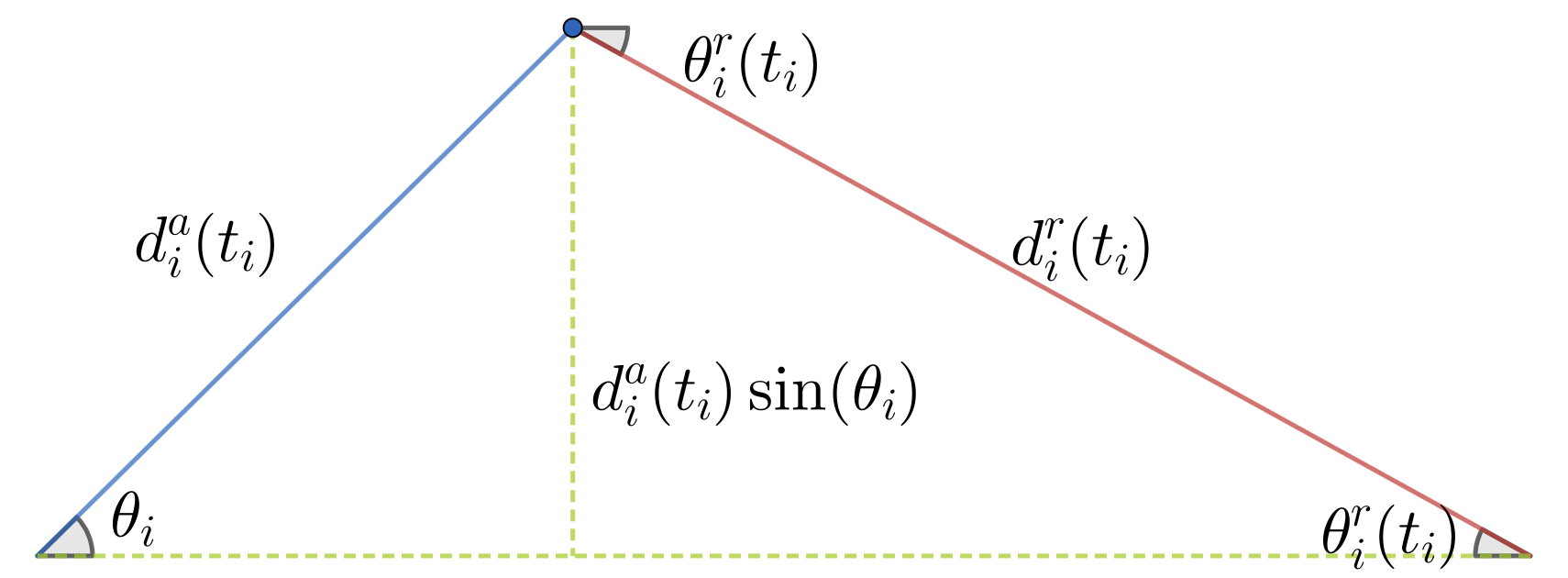}
    \caption{Calculation of the recovery angle $\theta^r_i$. The segment in blue corresponds to the avoidance stage while the segment in red corresponds to the trajectory recovery stage. The segments in green are the projection of trigonometric functions in a right triangle.}
    \label{fig:rec3}
\end{figure}

In this study, we are considering a simple trajectory recovery model in which deviated aircraft perform a second manoeuvre in order to change their trajectory towards their target point. Let $(\check{x}_i,\check{y}_i)$ be the coordinate of the target point, those second manoeuvres are opposing to the deviations in the avoidance stage. The speed component $q^r$ is simply defined as the aircraft are returning their velocity profile to its nominal value, i.e. $q^r_i = 1$. For the heading changes, the deviation angle is based on the time each aircraft moves from its avoidance trajectory towards its target point. For a given aircraft $i \in \mathcal{A}$, its recovery trajectory is defined as : $p_i(t) = [\widehat{x}_i + q_i\hat{v}_i\cos (\widehat{\theta}_{i} + \theta_i)t_i + \hat{v}_i\cos (\widehat{\theta}_{i})t, \widehat{y}_i + q_i\hat{v}_i\sin (\widehat{\theta}_{i} + \theta_i)t_i + \hat{v}_i\sin (\widehat{\theta}_{i})t]^{T}$, where $t_i$ corresponds to the recovery time, i.e., the value of time in which each aircraft change from its avoidance trajectory to its trajectory recovery (see Figure \ref{fig:rec3}). Hence, the deviation angle $\theta^r_i$ can be calculated as:

\begin{equation}
\theta^r_i(t_i) = \arcsin\Big(\frac{d^a_i(t_i)\sin(\theta_{i})}{d^r_i(t_i)}\Big),
\end{equation}

where the $d^a_i(t_i)$ (see Eq. \eqref{da}) corresponds to the distance flown during the conflict avoidance stage (from the initial position until the aircraft reach $t_i$ and changes to trajectory recovery) and $d^r_i(t_i)$ (see Eq. \eqref{dr} to the distance flown during the recovery stage (from the time $t_i$ where aircraft starts its trajectory until it reaches its destination point $(\check{x}_i,\check{y}_i)$) (see Figure \ref{fig:rec3}). If the aircraft did not change its heading angle, the recovery angle is not calculated.

\begin{equation}
    \label{da}
    d^a_i(t_i) = \sqrt{(\widehat{x}_i - x(t_i))^2 + (\widehat{y}_i - y(t_i))^2}, \forall i \in \mathcal{A},
\end{equation}

\begin{equation}
    \label{dr}
    d^r_i(t_i) = \sqrt{(x(t_i) - \check{x}_i)^2 + (y(t_i) - \check{y}_i)^2}, \forall i \in \mathcal{A}.
\end{equation}

In the avoidance stage, the core idea is to compare the distance between each pair of aircraft and derive time-independent expressions that can be escalated to all instances and applied simultaneously. However, the same cannot be achieved when avoidance and recovery are calculated simultaneously. This is because the aircraft recovery trajectories are a function of the time when aircraft perform their recovery manoeuvre ($t_i$). It is evident that this expression is nonlinear concerning speed and heading variables. This solidifies that, by using mathematical programming, it is quite complex to reproduce the avoidance model to solve the recovery stage. 

\begin{figure}[!h]	
	\centering	
		\subfloat[CP-5: circle problem with 5 aircraft]{\includegraphics[width=0.5\textwidth]{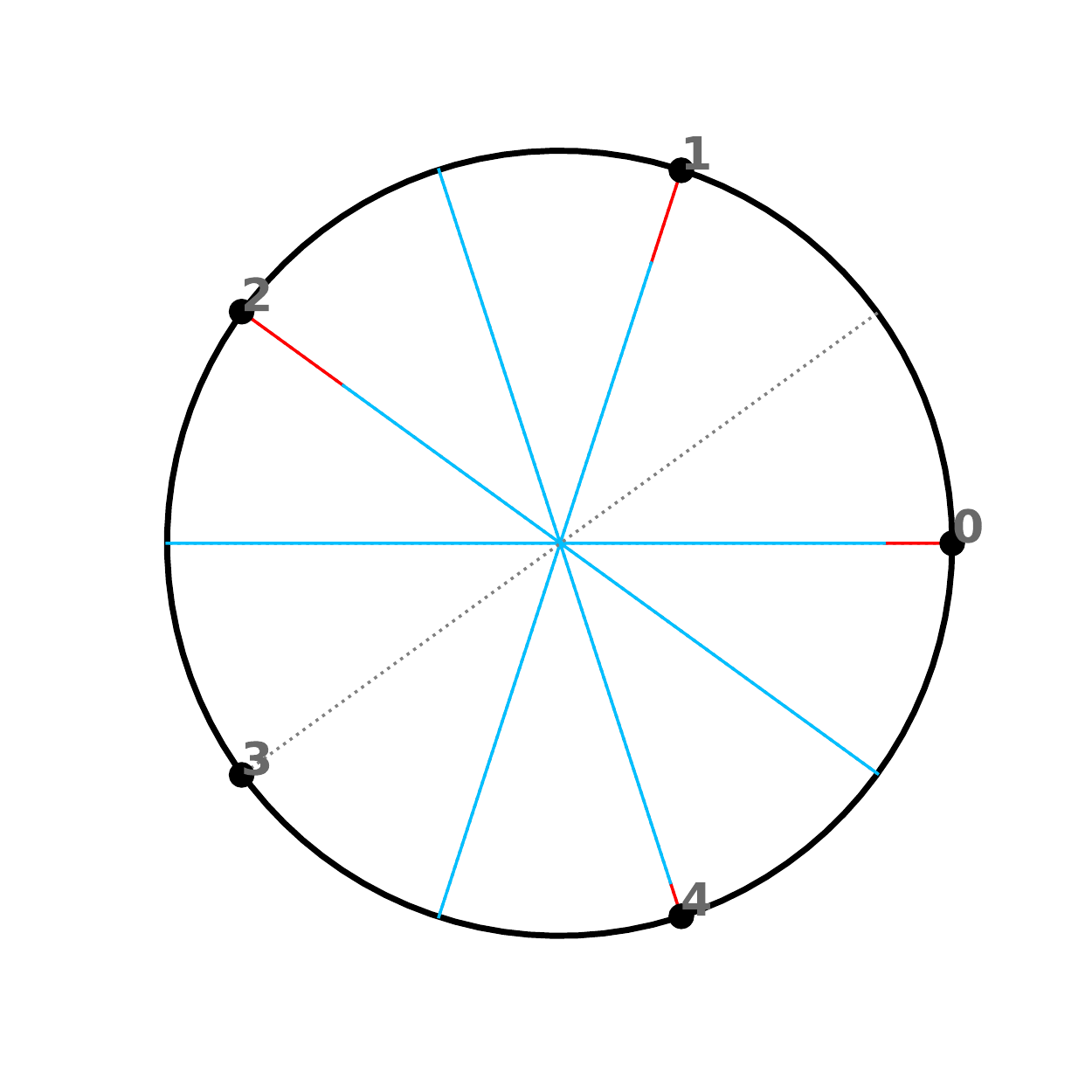}\label{fig:cp-5}} %\hspace*{2cm}
		\subfloat[RCP-10: random circle problem with 10 aircraft]{\includegraphics[width=0.5\textwidth]{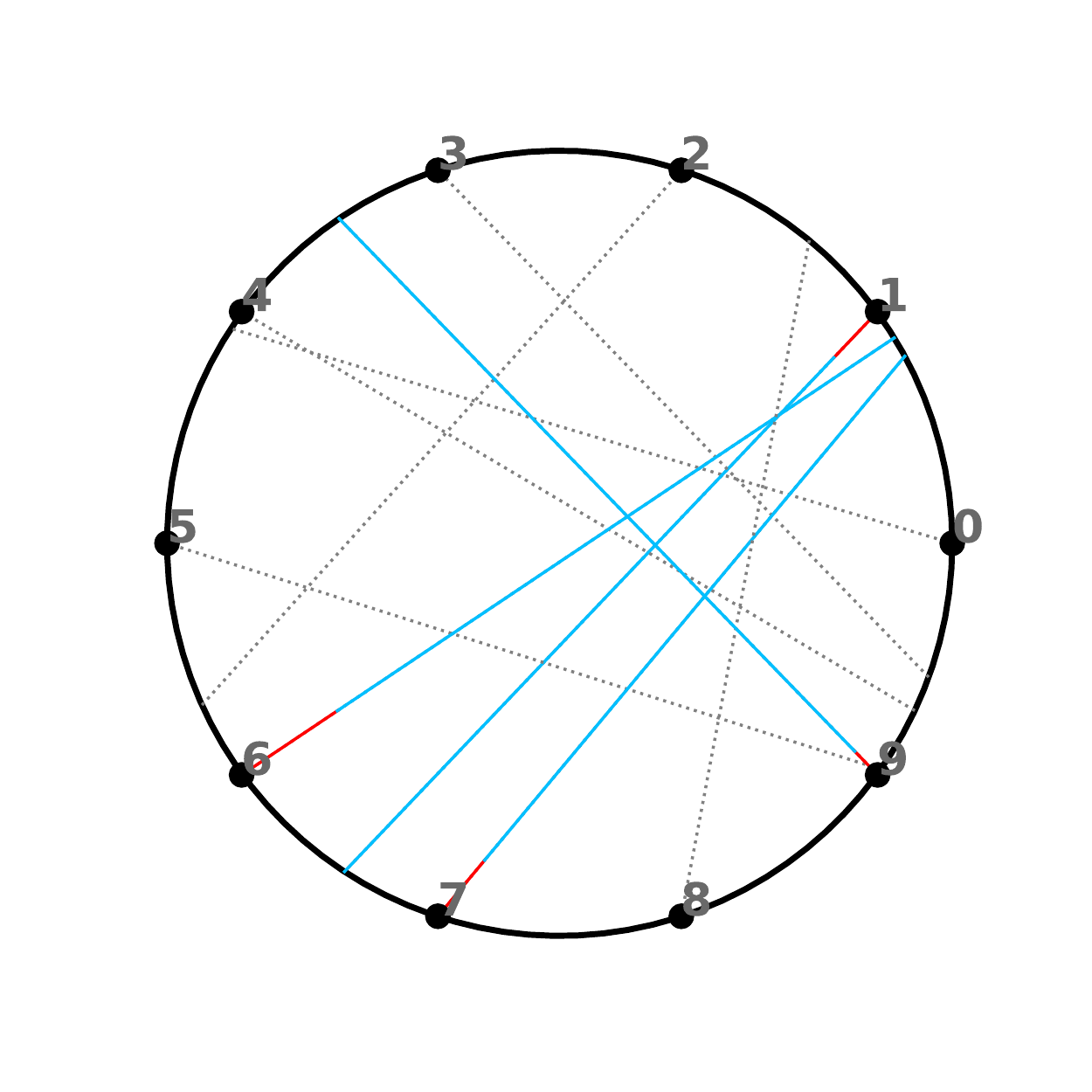} \label{fig:CP-7}}\\
		\subfloat[CP-15: circle problem with 15 aircraft]{\includegraphics[width=0.5\textwidth]{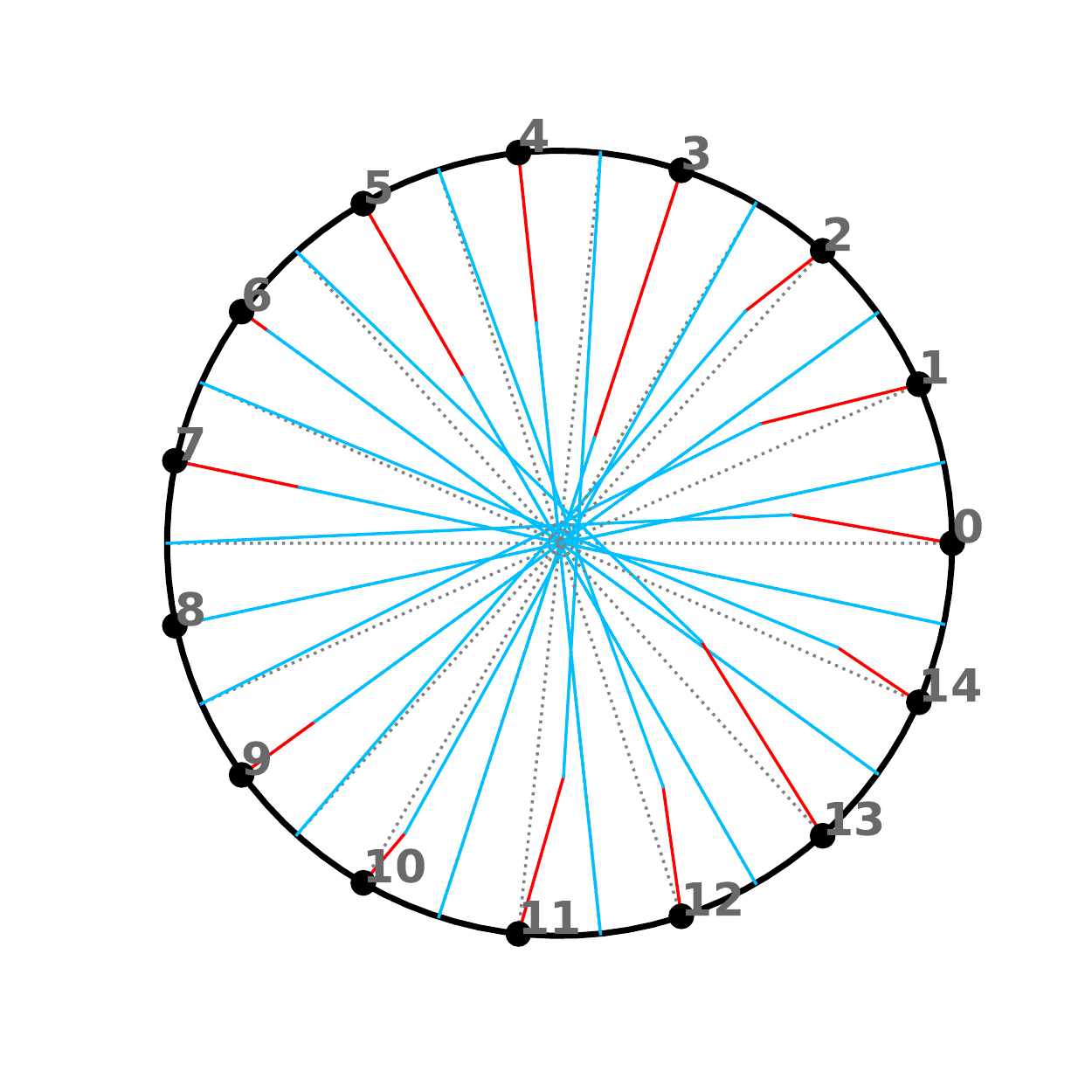}\label{fig:cp-15}} %\hspace*{2cm}
		\subfloat[RCP-20: random circle problem with 20 aircraft]{\includegraphics[width=0.5\textwidth]{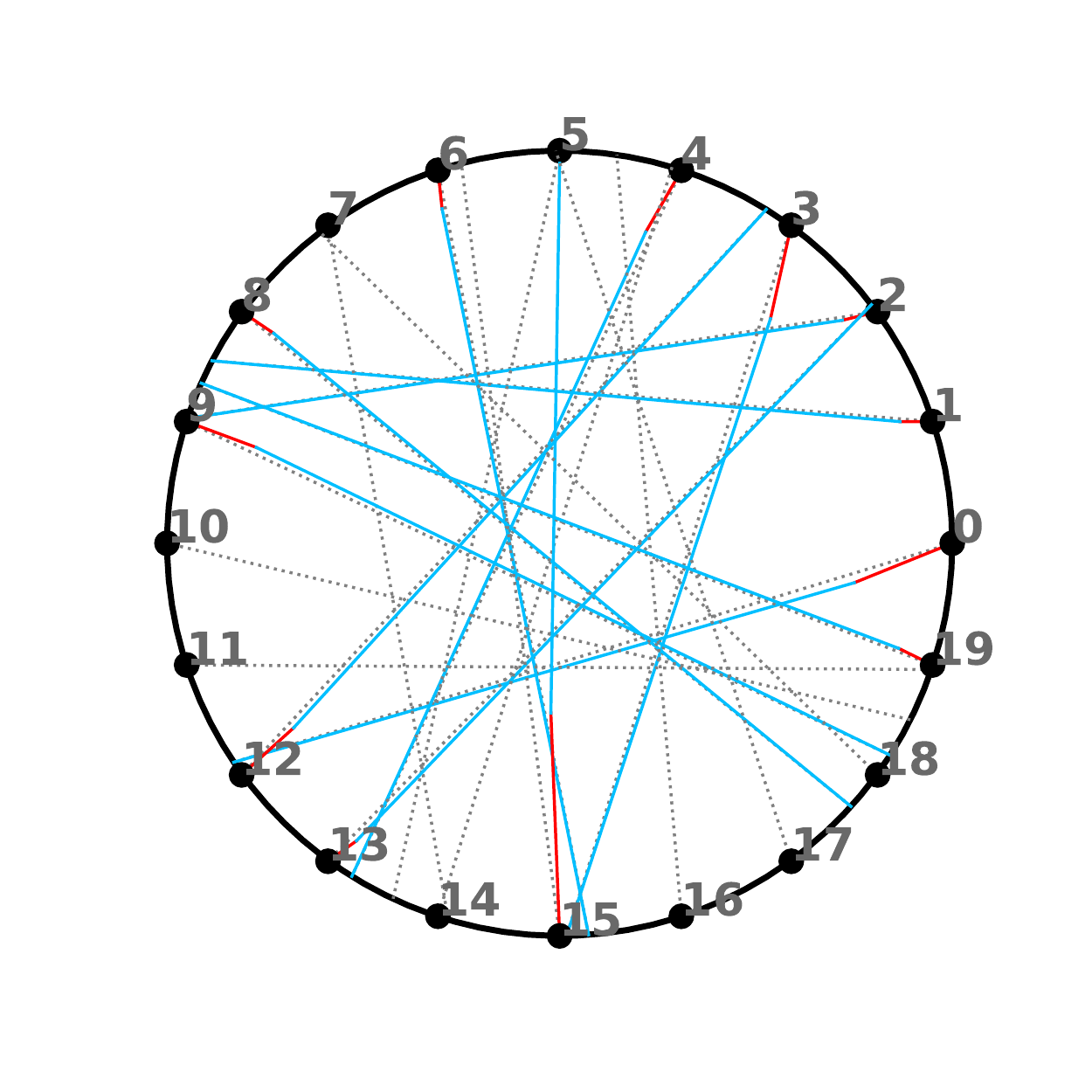} \label{fig:RCP-20-1}}\\
		\hfil
	\caption{Example of benchmarking instances for the Circle Problem (CP) using the naive approach in \cite{dias2020two}}.
	\label{fig:naive}
\end{figure}

Given this context, it is expected that some simplification level is required to make any formulation using mathematical programming viable. One of these strategies were implemented by \cite{dias2020two}. In their paper, a naive approach for trajectory recovery was implemented, and it is composed of a two-stage algorithm where avoidance and recovery are solved sequentially. In this formulation, the manoeuvres determined during the action stage are compensated in the recovery stage. In order to incorporate the cost of avoidance in the recovery stage, the total deviation on avoidance is passed on to the recovery stage. It determines how costly the deviation is in terms of speed and angle applied. In this way, the trajectory recovery can compensate for the cost during the avoidance stage. Another characteristic of this formulation is that the recovery time is a discrete variable. This reduces the feasible region and makes the options for trajectory recovery a limited, finite set. Based on those characteristics, this formulation can solve small to more significant instances (up to 30 aircraft) in a reasonable amount of time. More details can be found in \cite{dias2020two}. An example of the solution obtained via this approach can be seen in Figure \ref{fig:naive}.

The main drawback of this formulation is the lack of anticipation of the recovery costs at the avoidance stage. As both stages are solved separately, and each stage is solved only once, the solution obtained after solving both stages is myopic. This is not a significant issue for small instances due to the small amount of aircraft, but it becomes concerning in medium to large instances. As the number of aircraft and conflicts increase, the naive approach of \cite{dias2021disjunctive} is expected to lead to weak solutions that do not anticipate recovery costs in the design of aircraft trajectories.

In order to explore the situation where the avoidance takes into account the recovery costs, we are proposing an alternative version of such an algorithm where the total cost is optimized throughout an iterative algorithm by projecting the cost of recovery operations into the avoidance stage to obtain non-trivial solutions. In the following section, we explain this algorithm in detail.

\begin{figure}
	\centering
	\includegraphics[width=1.0\linewidth,width=0.7\linewidth]{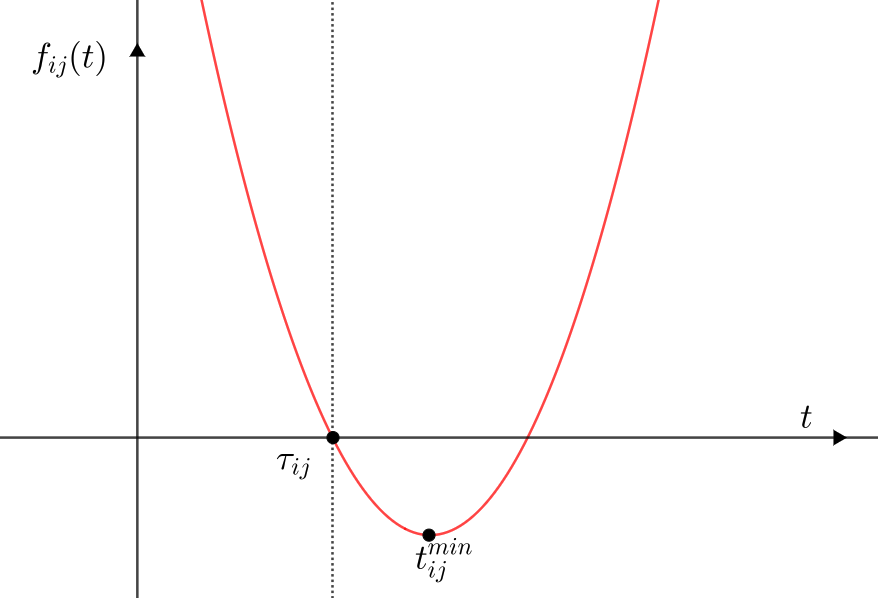}
	\caption{Illustration of $f_{ij}(t)$ for a configuration with $g_{ij} < 0$ and $\tm > 0$. $\tau_{ij}$ represents the start time of the conflict.}
	\label{fig:conflictandturn}
\end{figure}

Each aircraft needs to perform opposing manoeuvres for trajectory recovery to cancel the deviation applied during avoidance. Similar to the avoidance model, our goal is to guarantee that all pairs of aircraft are separated throughout the recovery stage. Since the separation condition in Eq. \eqref{eq:zd} is based on linear motion, we need to distinguish the trajectory stage of each aircraft $i \in \mathcal{A}$, i.e. before and after its recovery time $t_i$. 
We denote $A_i$ the avoidance trajectory of aircraft $i$ and $R_i$ its recovery trajectory. 
Given a pair $(i,j)$ of aircraft, we need to ensure that aircraft are separated during all pairwise trajectory stages, denoted $A_iA_j$, $A_iR_j$, $R_iA_j$ and $R_iR_j$. 
Observe that separation for the stage $A_iA_j$ is already ensured by the solution of Model \ref{mod_disj_rec}. If aircraft $i$ and $j$ were to recover at the same time period, then aircraft will transition from $A_iA_j$ to $R_iR_j$ directly. Otherwise, if $i$ (resp. $j$) recovers before $j$ (resp. $i$), then $A_iA_j$ will transition to $R_iA_j$ (resp. $A_iR_j$) before transitioning to $R_iR_j$. 
To avoid trigonometric forms and to obtain a tractable formulation, we discretize aircraft recovery times. Let $\mathcal{T}$ be the set of time periods available for recovery, we require:  

\begin{equation}
    t_i \in \{0,1\epsilon, 2\epsilon, \ldots, |\mathcal{T}|\epsilon\},
\end{equation}

where $\epsilon$ is the length of time periods. Abusing notation, we redefine the separation condition expressed in Eq. \eqref{eq:sepconditions} as: $g_{ij}(m,n) \geq 0$ and $\tm(m,n) \leq 0$ where the pair $(m,n)$ indicates the time period indices of recovery times $t_i$ and $t_j$, respectively. Let $\Omega_{X_iX_j}$ be the set of conflict-free pairs of recovery times for aircraft $i,j \in \mathcal{A}$ where $X_i$ represents the state of the trajectory of aircraft $i$, \emph{i.e.} $A_i$ or $R_i$; and $X_j$ represents the state of the trajectory of aircraft $j$, \emph{ i.e.} $A_j$ or $R_j$. This set can be specified into three different sets corresponding to the three different states during the recovery stage. The set $\Omega_{R_iR_j}$ is defined as:
\begin{align}
\Omega_{R_iR_j} &= \{(m,n)\in \mathcal{T}^2 :  g_{R_iR_j}(mn) \geq 0 \vee \underline{t}_{R_iR_j}(m,n)\leq 0\}.
\end{align}
For the states $A_iR_j$ and $R_iA_j$ an extra condition is required. Consider the state $A_iR_j$: if the lines of motion corresponding to trajectories $A_i$ and $R_j$ are in conflict but aircraft $i$ turns into recovery prior to the start of this conflict, then no conflict will occur. This illustrated in Figure \ref{fig:conflictandturn} where $g_{A_iR_j} < 0$ and $t_{A_iR_j} > 0$. Let $\tau_{A_iR_j}(t_j)$ be the smallest root of $g_{A_iR_j} = 0$ if $j$ recovers at time $t_j$. If aircraft $i$ recovers prior to $\tau_{A_iR_j}(t_j)$, \textit{i.e.} $t_i \leq \tau_{A_iR_j}(t_j)$, then the conflict will be avoided. Accordingly, we define: 
\begin{subequations}
\begin{align}
\Omega_{A_iR_j} &= \{(m,n) \in \mathcal{T}^2 :  g_{A_iR_j}(n) \geq 0  \vee t^{min}_{A_iR_j}(n)\leq 0  \vee m \leq \tau_{A_iR_j}(n)\}, \\
\Omega_{R_iA_j} &= \{(m,n) \in \mathcal{T}^2 :  g_{R_iA_j}(m) \geq 0 \vee t^{min}_{R_iA_j}(m)\leq 0  \vee  n \leq \tau_{R_iA_j}(m)\}.
\end{align}
\end{subequations}

Let $\rho_{im}$ be a binary variable equal to 1 if aircraft $i \in \mathcal{A}$ recovers at time period $m \in \mathcal{T}$ and 0 otherwise. To track the states of aircraft pair $(i,j)$ which are activated, we introduce two binary variables $\alpha_{ij}$ and $\beta_{ij}$. Those variables are used to identify whether $t_i < t_j$ ($\alpha_{ij} = 1$) which activates state $R_iA_j$, or if $t_i > t_j$ ($\beta_{ij} = 1$) which activates state $A_iR_j$. Variables $\alpha_{ij}$ and $\beta_{ij}$ are defined via the constraints:
\begin{subequations}\label{cons1}
	\begin{align}
	& \alpha_{ij} \geq \frac{1}{|\mathcal{T}|}\Bigg(\sum\limits_{m \in \mathcal{T}}m\rho_{im} - \sum\limits_{n \in \mathcal{T}}n\rho_{jn}\Bigg), &&\forall (i,j) \in \mathcal{P}, \\
	& \beta_{ij} \geq \frac{1}{|\mathcal{T}|} \Bigg(\sum\limits_{n \in \mathcal{T}}n\rho_{jn} - \sum\limits_{m \in \mathcal{T}}m\rho_{im} \Bigg), &&\forall (i,j) \in \mathcal{P} \\
	& \alpha_{ij} + \beta_{ij} \leq 1, &&\forall (i,j) \in \mathcal{P}.
	\end{align}
\end{subequations}

We use the following constraints to exclude conflicting trajectories from the solution. Observe that states $A_iR_j$ and $R_iA_j$ are conditional on the recovery times of $t_i$ and $t_j$ and thus the corresponding constraints are only active if $i$ and $j$ do not recover at the time period.
\begin{subequations}\label{cons2}
	\begin{align}
	& \rho_{im} + \rho_{jn} \leq 2 - \beta_{ij} &&\forall (i,j) \in \mathcal{P}, (m,n) \in \Omega_{A_iR_j},\\
	& \rho_{im} + \rho_{jn} \leq 2 - \alpha_{ij} &&\forall (i,j) \in \mathcal{P}, (m,n) \in \Omega_{R_iA_j}, \\
	& \rho_{im} + \rho_{jn} \leq 1  &&\forall (i,j) \in \mathcal{P}, (m,n) \in \Omega_{R_iR_j}.
	\end{align}
\end{subequations}
Aircraft are assigned a recovery time via the constraint:
\begin{align}
& \sum_{m  \in \mathcal{T}} \rho_{im} = 1 \quad \quad \forall i \in \mathcal{A},
\label{cons3}
\end{align}

The second stage aims to identify the optimal time for aircraft to recover towards their target position. The avoidance cost function of aircraft $i$, denoted $a_i$, is defined as:

\begin{equation}
    \label{fobj}
    a_i = (1-w)(1-q_i)^2 + w\theta_i^2 + \lambda_f f_i \quad \forall i \in \mathcal{A},
\end{equation}

For the recovery stage, we are proposing to minimize the total weighted recovery time, which accounts for manoeuvres applied in the avoidance stage as well as aircraft recovery times.

\begin{equation}
\sum_{i \in \mathcal{A}} \sum_{m  \in \mathcal{T}} a_i\rho_{im}t_m^2,
\end{equation}

this expression comes from the approximated area created by the deviation during the avoidance stage, and it is quadratic in $t_i^2$. The trajectory recovery formulation is summarised in Model \ref{mod:onestepdetailed} which is a MIQP (Mixed Integer Quadratic Programming). 

\begin{model}[Trajectory Recovery]
	\label{mod:onestepdetailed}
	\allowdisplaybreaks
	\begin{align*}
	&\text{\emph{minimise}}\qquad \sum_{i \in \mathcal{A}} \sum_{m  \in \mathcal{T}} a_i\rho_{im}t_m^2\\
	&\text{\emph{subject to}} \nonumber \qquad \\
	& \alpha_{ij} \geq \frac{1}{|\mathcal{T}|}\Bigg(\sum\limits_{m \in \mathcal{T}}m\rho_{im} - \sum\limits_{n \in \mathcal{T}}n\rho_{jn}\Bigg) &&\forall (i,j) \in \mathcal{P}, \\
	& \beta_{ij} \geq \frac{1}{|\mathcal{T}|} \Bigg(\sum\limits_{n \in \mathcal{T}}n\rho_{jn} - \sum\limits_{m \in \mathcal{T}}m\rho_{im} \Bigg) &&\forall (i,j) \in \mathcal{P}, \\
	& \alpha_{ij} + \beta_{ij} \leq 1 &&\forall (i,j) \in \mathcal{P},\\
	& \rho_{im} + \rho_{jn} \leq 2 - \beta_{ij} &&\forall (i,j) \in \mathcal{P}, (m,n) \in \Omega_{A_iR_j},\\
	& \rho_{im} + \rho_{jn} \leq 2 - \alpha_{ij} &&\forall (i,j) \in \mathcal{P}, (m,n) \in \Omega_{R_iA_j}, \\
	& \rho_{im} + \rho_{jn} \leq 1  &&\forall (i,j) \in \mathcal{P}, (m,n) \in \Omega_{R_iR_j}, \\
	& \sum_{m  \in \mathcal{T}} \rho_{im} = 1 &&  \forall i \in \mathcal{A}, \\
	& \rho_{im} \in \{0,1\}  && \forall i \in \mathcal{A}, m \in \mathcal{T}, \\
	& \alpha_{ij},\beta_{ij} \in \{0,1\}  && \forall (i,j) \in \mathcal{P}.
	\end{align*}
\end{model}

\subsection{Penalty-based Conflict Resolution and Trajectory Recovery Algorithm}

The main idea behind the proposed penalty-based approach is to capture the recovery cost during the avoidance stage. Therefore, its goal is to influence the behaviour of the avoidance stage to anticipate the cost of trajectory recovery and attempt to construct an efficient trajectory across both stages. This can be achieved by using the solution of the recovery stage as a preemptive cost in the avoidance stage. In this case, the algorithm aims to find a trade-off between the deviation costs incurred during the avoidance and recovery stages. Let $TC_i$ be the total cost per aircraft defined as the combined cost of avoidance and recovery. 

\begin{equation}
    TC_i = (1-w)(1-q_i)^2 + w\theta_i^2 + \lambda_f f_i + \lambda_t t_i^2,\qquad \forall i \in \mathcal{A},
    \label{totalcost}
\end{equation}

where $\lambda_t$ is the weight for the recovery time component. 

Therefore, the objective function in the avoidance stage is modified to account for the anticipated cost of trajectory recovery. Let $r_i$ the recovery cost of each aircraft $i \in \mathcal{A}$, where the avoidance cost is based on the deviation in the avoidance stage ($q_i^\star$ and $\theta_i^\star$ being the optimal values for speed changes and heading angles).

\begin{equation}
    \label{eq:recovery}
    r_i = (1-w)(1-q_i^\star)^2 + w(\theta_i^\star)^2 + \lambda_t t_i^2, \forall i \in \mathcal{A},
\end{equation}

By adding this expression into the avoidance stage, we aim to account for the impact of recovery costs within the avoidance stage. This approach aims to penalize the control of aircraft in the avoidance stage proportionally to their anticipated total trajectory deviation. Thus, the objective function in the avoidance stage can be rewritten as:

\begin{equation}
    \min \sum_{i \in \mathcal{A}} w\theta_i^2 + (1 - w)(1 - q_i)^2 + f_i(\lambda_f + {r}_i)
    \label{total}
\end{equation}

The proposed iterative approach is summarized in Algorithm \eqref{algo:2d_rec}. In each iteration, the avoidance stage (as described in Model \eqref{mod_disj_rec}) is solved using the updated objective function described in Eq. \eqref{total}. 

Based on the solutions obtained from solving Model \eqref{mod_disj_rec}, the sets $\Omega_{A_iR_j},\Omega_{AjR_i}$ and $\Omega_{R_iR_j}$ are pre-calculated. The recovery stage is solved subsequently and the total cost and the cost variation $\deltatc$ (based on Eqs. \eqref{totalcost} and \eqref{deltatc}) are calculated. The stop criteria is based on the overall total cost between two consecutive iterations, i.e.:

\begin{equation}
    \label{deltatc}
    \Delta TC = \sum_{i \in \mathcal{A}} TC_i^{n} - TC_i^{n-1},
\end{equation}

where $n$ corresponds to the number of the iteration. If the value of $\Delta TC$ is below a predefined threshold, the algorithm converges and stops; otherwise, the recovery cost is updated and the algorithms proceeds. In the first iteration, we assume $r_i = 0$; for the consecutive iterations, the anticipated recovery cost $r_i$ is calculated using the solution obtained in the previous iteration.

\begin{algorithm}
\caption{Solution algorithm for the Conflict Resolution Problem with Trajectory Recovery}\label{algo:2d_rec}
\begin{algorithmic}[1]
\Require $\mathcal{A},\hat{\bm{\theta}}, \hat{\bm{v}}, \qlb, \qub, \tlb, \tub, \epsilon$
\Ensure $\bm{q}^\star, \bm{\theta}^\star, \bm{t}_i^\star, \bm{f}_i^\star$, LB, UB
\State LB, UB
\State LB $\gets 0$ 
\State UB $\gets +\infty$
\State gap $\gets 0$ 
\State \textsc{converged} $\gets$ \texttt{False}
\While{\textsc{converged} $=$ \emph{\texttt{False}}}
        \State $\bm{q}, \bm{\theta}, \bm{z},\bm{f}$, LB $\gets$ Solve Stage 1 - Avoidance using Model \eqref{mod_disj_rec}
        \If{\emph{Infeasible}}
            \State \Return \textsc{Infeasible}
        \EndIf
        \State Calculate sets $\Omega_{A_iR_j},\Omega_{AjR_i}$ and $\Omega_{R_iR_j}$ according to Eqs. \eqref{cons2}
        \State Calculate set $\Po$ 
        \State Add cuts in Eq. \eqref{fcons} 
	    \State $\bm{t}$ $\gets$ Solve Stage 2 - Recovery using Model \eqref{mod:onestepdetailed}
        \If{\emph{Infeasible}}
            \State \Return \textsc{Infeasible}
        \Else
        \State $\bm{t}^\star \gets \bm{t}$ 
    	\State $\bm{q}^\star \gets \bm{q}$
    	\State $\bm{\theta}^\star \gets \bm{\theta}$
    	\State $\bm{t}^\star \gets \bm{t}$
    	\State Calculate $\Delta TC$
        \EndIf
        \If{$\Delta TC \leq \epsilon$}
        \State \textsc{converged} $\gets$ \texttt{True}
        \Else
        \State $r_i \gets w\theta_i^2 + (1-w)(1-q_i)^2 + \lambda_tt_i^2$ 
        \EndIf
\EndWhile
\end{algorithmic}
\end{algorithm}

In order to solve Model \eqref{mod_disj_rec}, which is nonconvex and has nonlinear constraints that are notoriously challenging, we use the complex number formulation described by \cite{rey2017complex} and an adaptation of the solution method proposed by \cite{dias2021disjunctive}. 

\section{Numerical Results}
\label{num_contin}

The experimental framework used to test the proposed mixed-integer formulation for the trajectory recovery using continuous heading is introduced in Section \ref{exp_cont}. Then a detailed analysis of four groups of instances is presented in Section \ref{ilus_cont}. The computational performance of the proposed approaches is thoroughly explored in Section \ref{perf_cont}, respectively. 

\subsection{Experimental Framework}
\label{exp_cont}

\begin{figure}[t]	
	\centering	
		\subfloat[CP-7: circle problem with 7 aircraft]{\includegraphics[width=0.45\textwidth]{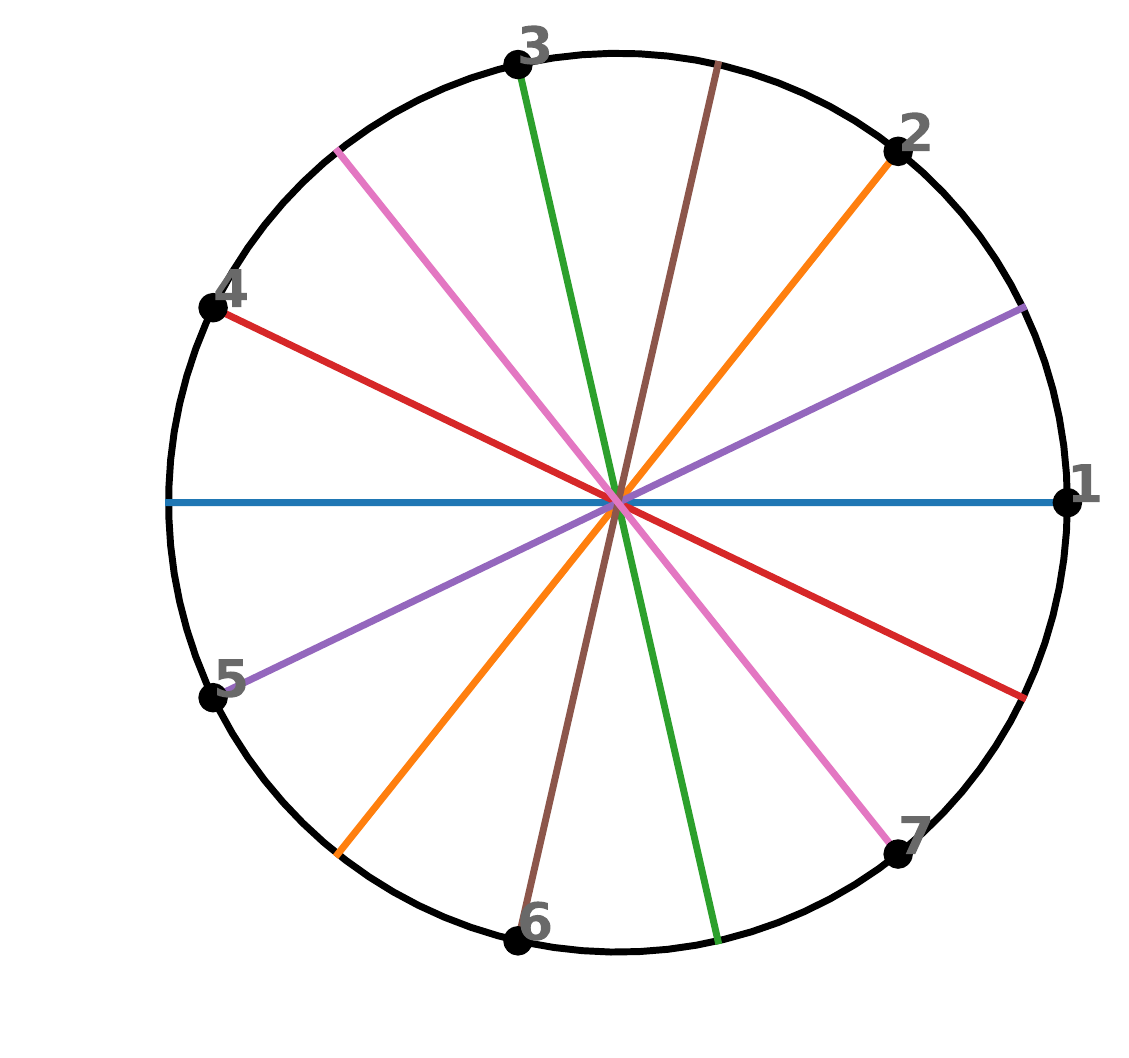}\label{fig:cp}} %\hspace*{2cm}
		\hfil
		\subfloat[RCP-10: random circle problem with 10 aircraft]{\includegraphics[width=0.45\textwidth]{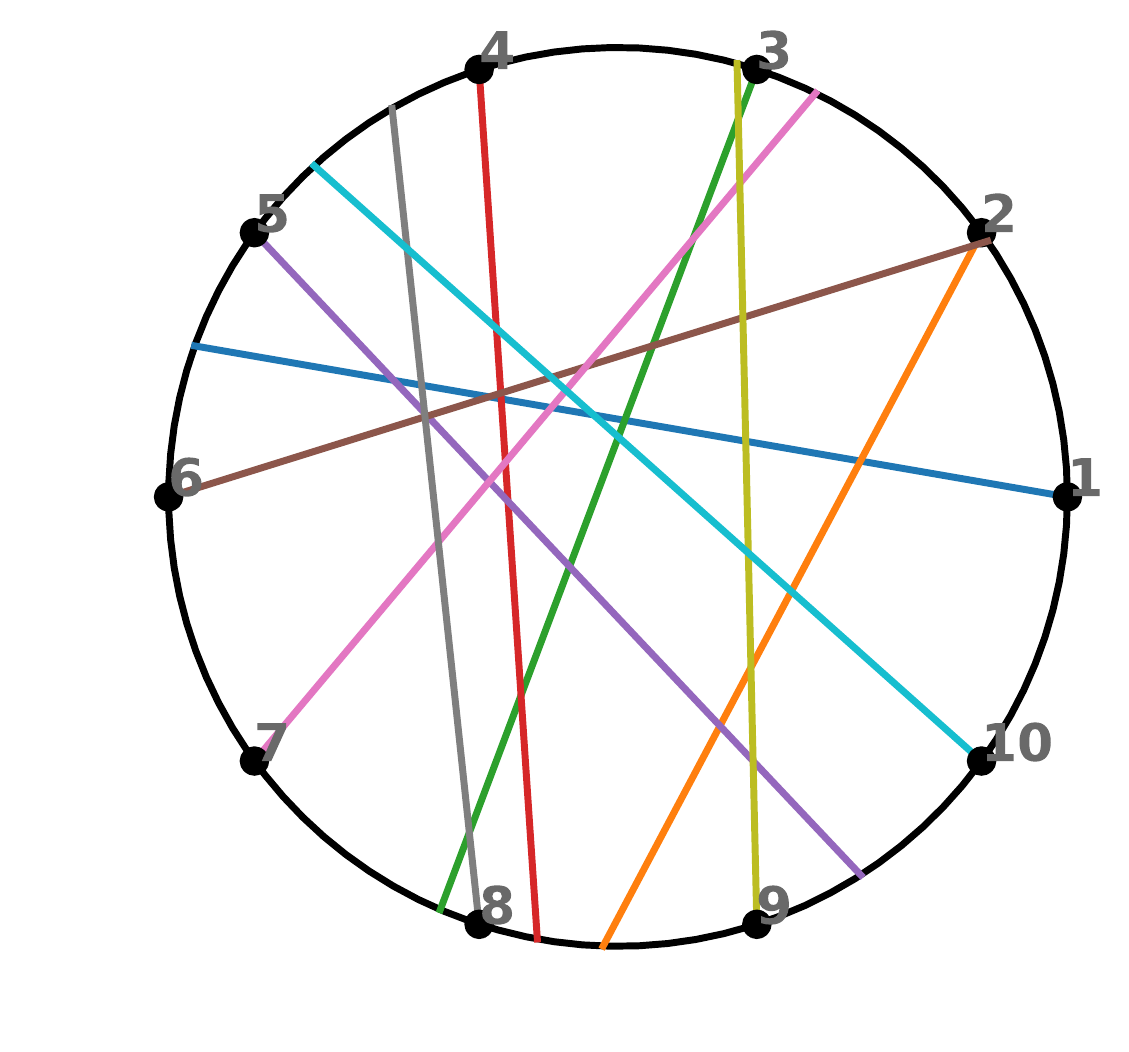} \label{fig:rcp10-5}}
	\caption{Example ofbenchmarking instances for the Circle Problem (CP) and Random Circle Problem (RCP).}
	\label{case1}
\end{figure}

We test the performance of the proposed mixed-integer formulations and algorithm using two benchmark problems from the literature: the Circle Problem (CP) and the Random Circle Problem (RCP). The two types of benchmarking instances are illustrated in Figure \ref{case1}. The CP consists of a set of aircraft uniformly positioned on the circle heading towards its centre. Aircraft speeds are assumed to be identical; hence the problem is highly symmetric (see Figure \ref{fig:cp}). The CP is notoriously difficult due to the geometry of initial aircraft configuration and has been widely used for benchmarking CD\&R algorithms in the literature \citep{durand2009ant,rey2015equity,cafieri2017mixed,cafieri2017maximizing,rey2017complex}. CP and RCP instances are named CP-N and RCP-N, respectively, where N is the total number of aircraft. For reproducibility purposes, all formulations and data used for testing are made available online at the public repository \small \url{https://github.com/acrp-lib/acrp-lib}.\normalsize 

For all tests, a speed regulation range based on the subliminal speed control $[-6\%,+3\%]$ is used, and the typical heading control range $[-30^\circ,+30^\circ]$ is also used. For the preference weight, we use as $w=0.5$ and $\lambda_f=1$ in the objective in Eq. \eqref{fobj}. This value was selected such that both heading and speed control terms were of a comparable order of magnitude, emphasizing penalizing heading control. For stage 2, a total of $|\mathcal{T}|=15$ time periods are used, with a step of $\epsilon=2$ minutes. To solve avoidance, the solution method used was an adaption of the algorithm implemented by \cite{dias2021disjunctive} \Dis, and for the recovery stage, the model is described in \ref{mod:onestepdetailed}. Both are solved with \Cplex Python API and a time limit of 5 minutes per solving and 15 minutes per instance. For the first iteration, $\rho_i^\star$ is assumed as $1$, and it is calculated after the recovery in the following iterations. 

The proposed approach is compared to the algorithms presented by \cite{dias2020two}. Those methods are respectively referred to as \Naive and \Greedy Recovery. The first method corresponds to a similar formulation of the method presented in this paper: it is the first iteration of Algorithm \ref{algo:2d_rec} where the avoidance and recovery stages are each solved only once. The second is a heuristic-based trajectory recovery procedure that iterates over all time steps and uses a priority list to decide which aircraft can be recovered at each time step. The priority list used is based on $r_i$ values \eqref{eq:recovery}. The algorithm first sorts aircraft accordingly and iterates overtime periods. At each time, the algorithm iterates over the sorted list of aircraft and check if each aircraft can be recovered at the current time. The process is repeated until no aircraft can recover at the current time. The proposed algorithm has a worst-case time complexity of $\mathcal{O}(|\mathcal{T}||\mathcal{A}|^3)$. Details of our implementation of both benchmarking algorithms can be found in \cite{dias2020two} and example of their solution can be found in Fig. \eqref{fig:sa}. Those models were slightly modified to be comparable to the algorithm presented in this paper. Specifically, the control variable $f_i$ is introduced in the avoidance stage. Since the solution methods presented in \cite{dias2020two} only solve each stage a single time, they are comparable to the first iteration of Algorithm \eqref{algo:2d_rec}.

\subsection{Illustration}
\label{ilus_cont}

\begin{figure}[htp]
    \centering
    \subfloat[CP-8]{\includegraphics[width=0.50\textwidth]{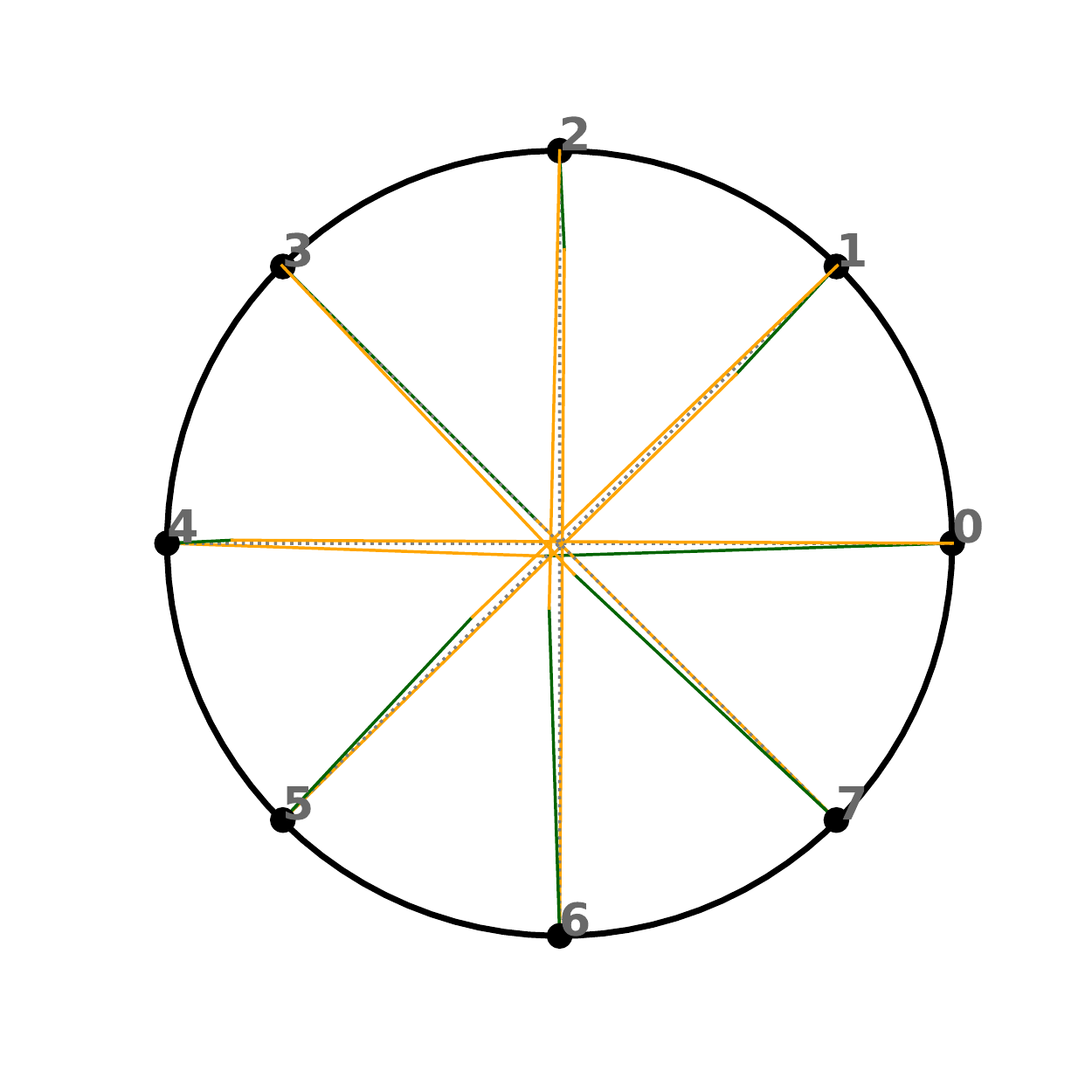}\label{cp7}}
    \hfil
    \subfloat[RCP-10]{\includegraphics[width=0.50\textwidth]{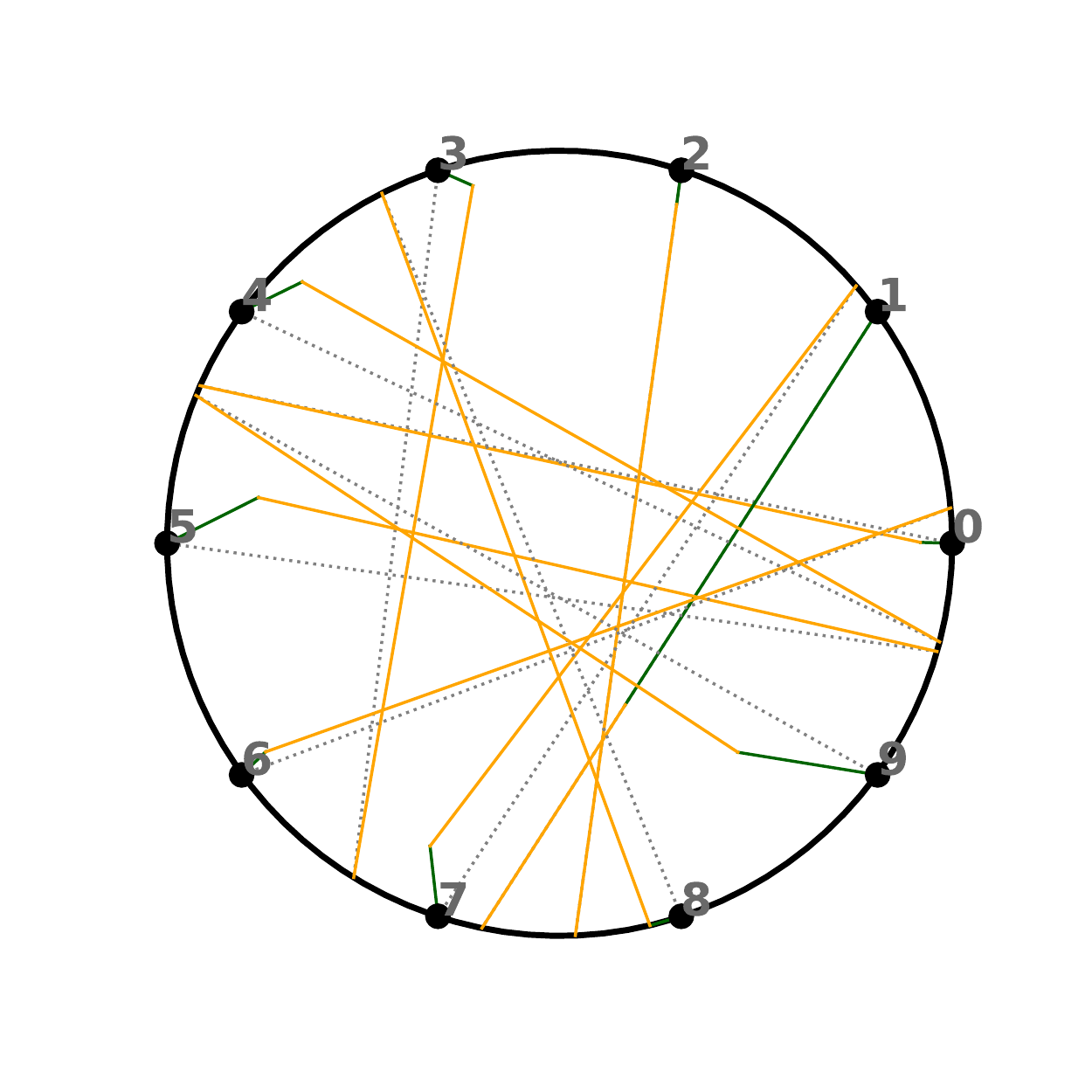}\label{cp10_6}}\\
    \subfloat[RCP-20]{\includegraphics[width=0.50\textwidth]{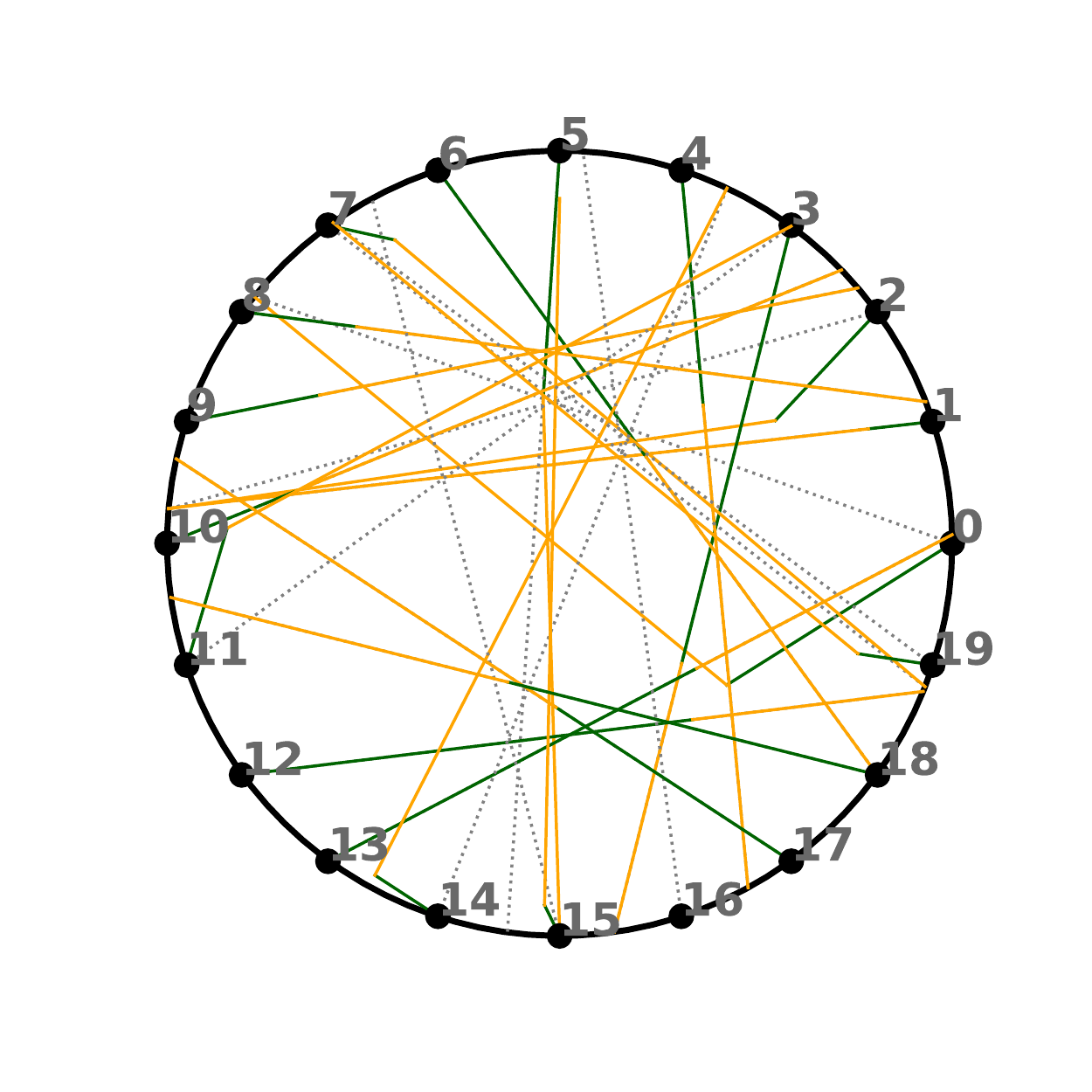}\label{gp10}}
    \hfil
	\subfloat[RCP-30]{\includegraphics[width=0.50\textwidth]{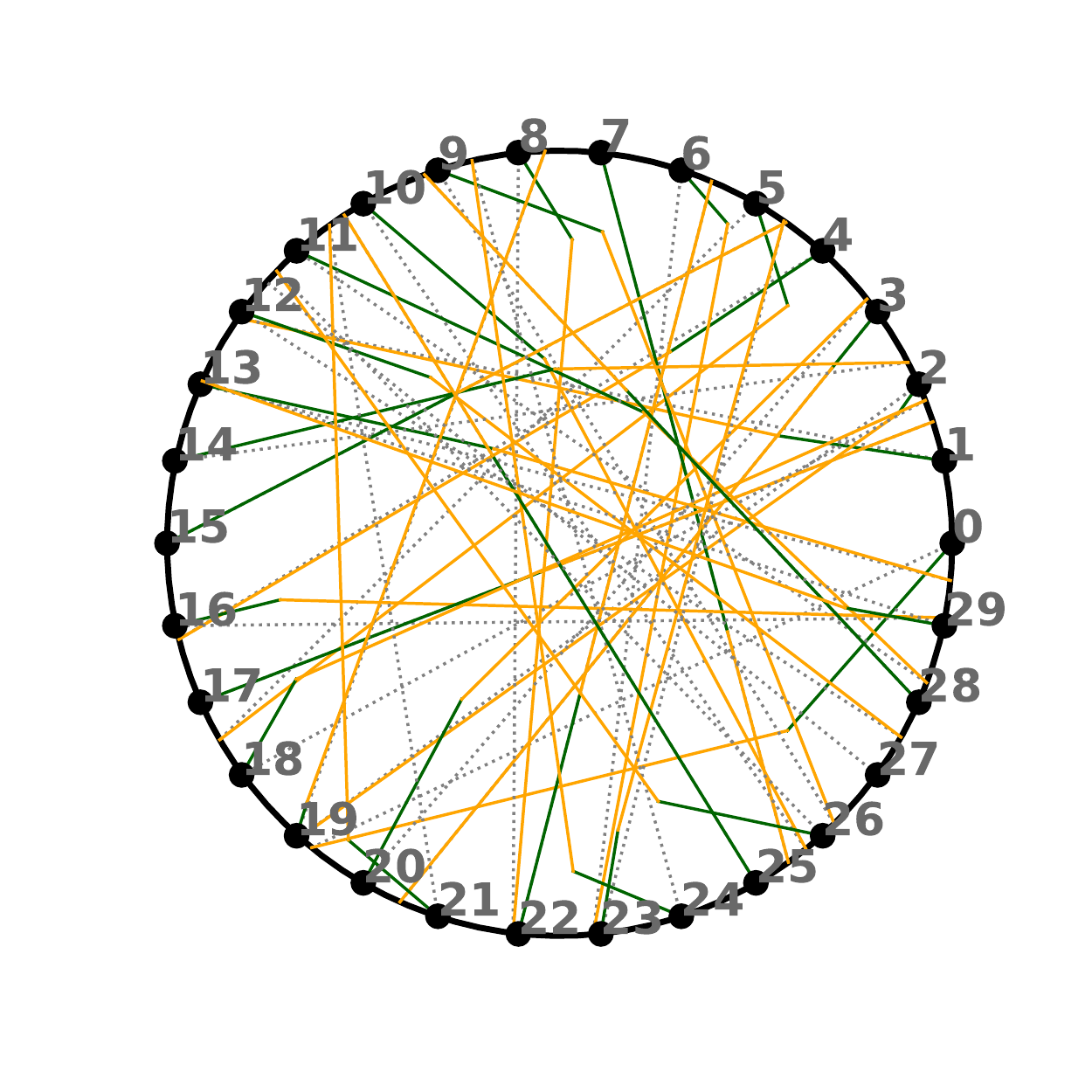}\label{rcp30}}\\
	\caption{Example of four instances used to benchmark our formulation. For all figures, lines in grey represents their nominal trajectory, the lines in orange represents the avoidance trajectory while the trajectories in green corresponds to their recovery. For those results are using special bounds for speed ([0.96,1.03]) and for heading angle([$\frac{-\pi}{3}$,$\frac{\pi}{3}$])}
	\label{fig:sa}
\end{figure}

To illustrate the proposed two-stage iterative algorithm, the optimal solution obtained is plotted for CP instances with five aircraft and RCP instances with 10, 20 and 30 aircraft. In Figure \eqref{fig:sa}, dashed grey lines represent initial aircraft trajectories, green lines represent the avoidance trajectory of stage 1, and orange lines represent recovery trajectories of stage 2 using Model \ref{mod:onestepdetailed}. In contrast, green lines represent the avoidance trajectory of stage 1, and the orange line represents the trajectory of stage 2 using Algorithm \ref{algo:2d_rec}.

This behaviour is even more accentuated in RCP-10-1, RCP2-20-2 and RCP-30-2. In the first, it is shown that only three aircraft are necessary to be altered with relatively minor deviations, and that would be sufficient for guaranteeing separation conditions instead of deviating 9 out of 10 aircraft. The solutions obtained via the continuous angle causes the slightest disturbance in the network, which is a favourable trade for those formulations. For RCP-20 and RCP-30, there are more aircraft to be altered, but overall, they show the same behaviour: more aircraft with higher avoidance cost and reduced recovery costs, with overall cost improved.

\subsection{Sensitivity Analysis}

To quantify the impact of the preference weights $\lambda_t$ and $\lambda_f$ in the proposed total cost function Eq. \eqref{totalcost}, we conduct numerical experiments on one instance of each of the two types of benchmarking instances for varying values of those parameters individually. This experiment focuses on the typical heading control range $[-30^\circ,+30^\circ]$. Our goal is to show that by varying those preference weight in $\lambda_{t} \in \ ]0,1[$ and $\lambda_{f} \in \ ]0,1[$ the decision-maker can control the desired level of trade-off between the recovery time and total deviation. Recall that in the total cost function \eqref{total}, $\lambda_f$ is the coefficient of $f_i$ which is minimal for $f_i = 0$; therefore, one can expect that increasing $\lambda_f$ will tend to penalize the number of aircraft controlled and therefore decrease the recovery time. On the other hand, increasing the value of $\lambda_t$ will tend to penalize the recovery time and, therefore, lead to potentially a greater total deviation of the aircraft.

\begin{figure}[!h]
    \centering
    \subfloat{\includegraphics[width=1.1\textwidth]{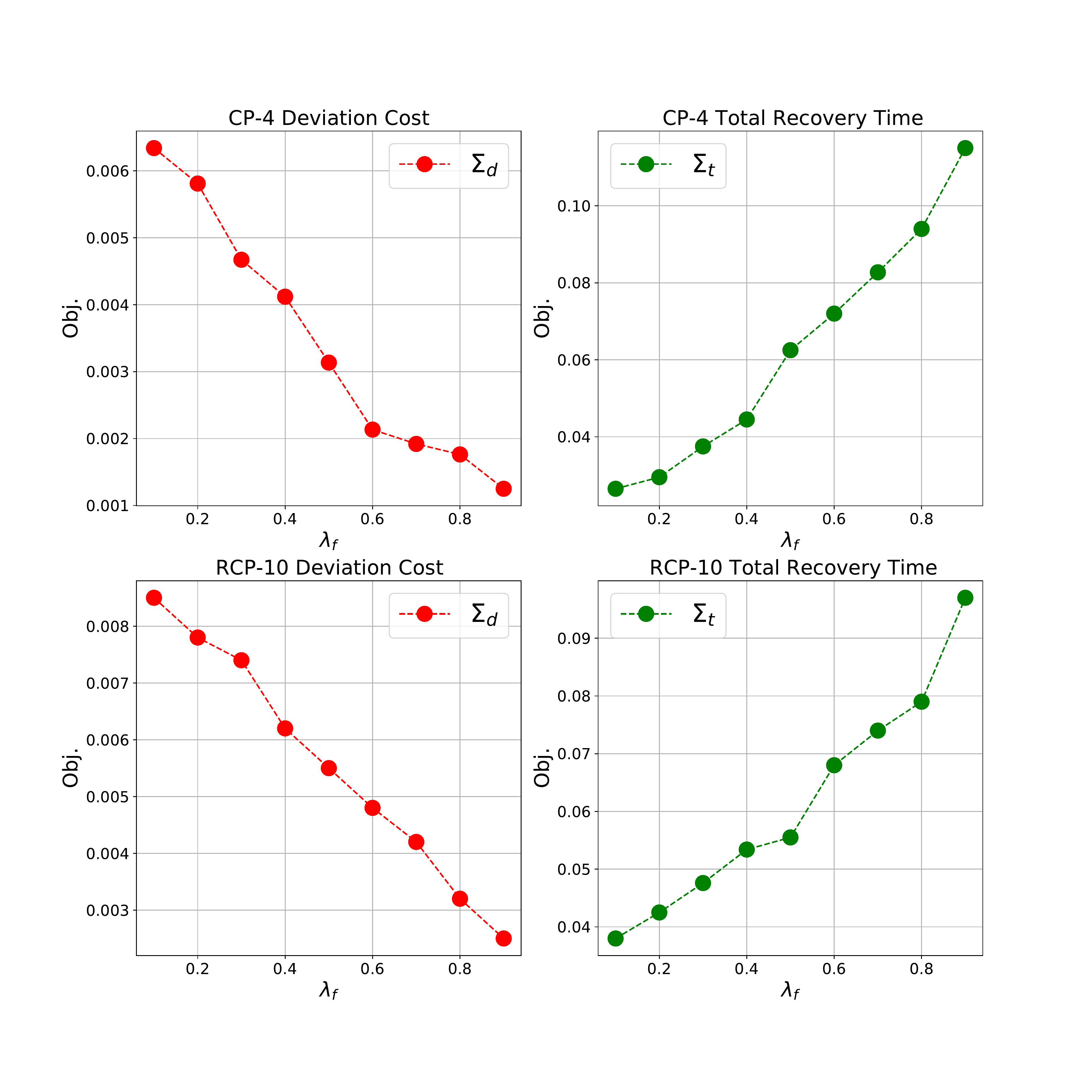}\label{scp7}}\\
	\caption{Sensitivity analysis on the preference weight $\lambda_f$ in the objective function Eq. \eqref{obj}. For all figures, $\Sigma_d$ represents the total speed and heading deviation defined as $\sum\limits_{i \in \mathcal{A}}(1-q_i)^2 + \theta_i^2$ (in red) and $\Sigma_t$ represents the total recovery time defined as $\sum\limits_{i \in \mathcal{A}} t_i^2$ (in green). All instances tested here are CP-4 and RCP-10.}
	\label{fig:sensitivityf}
\end{figure}

\begin{figure}[!h]
    \centering
    \subfloat{\includegraphics[width=1.1\textwidth]{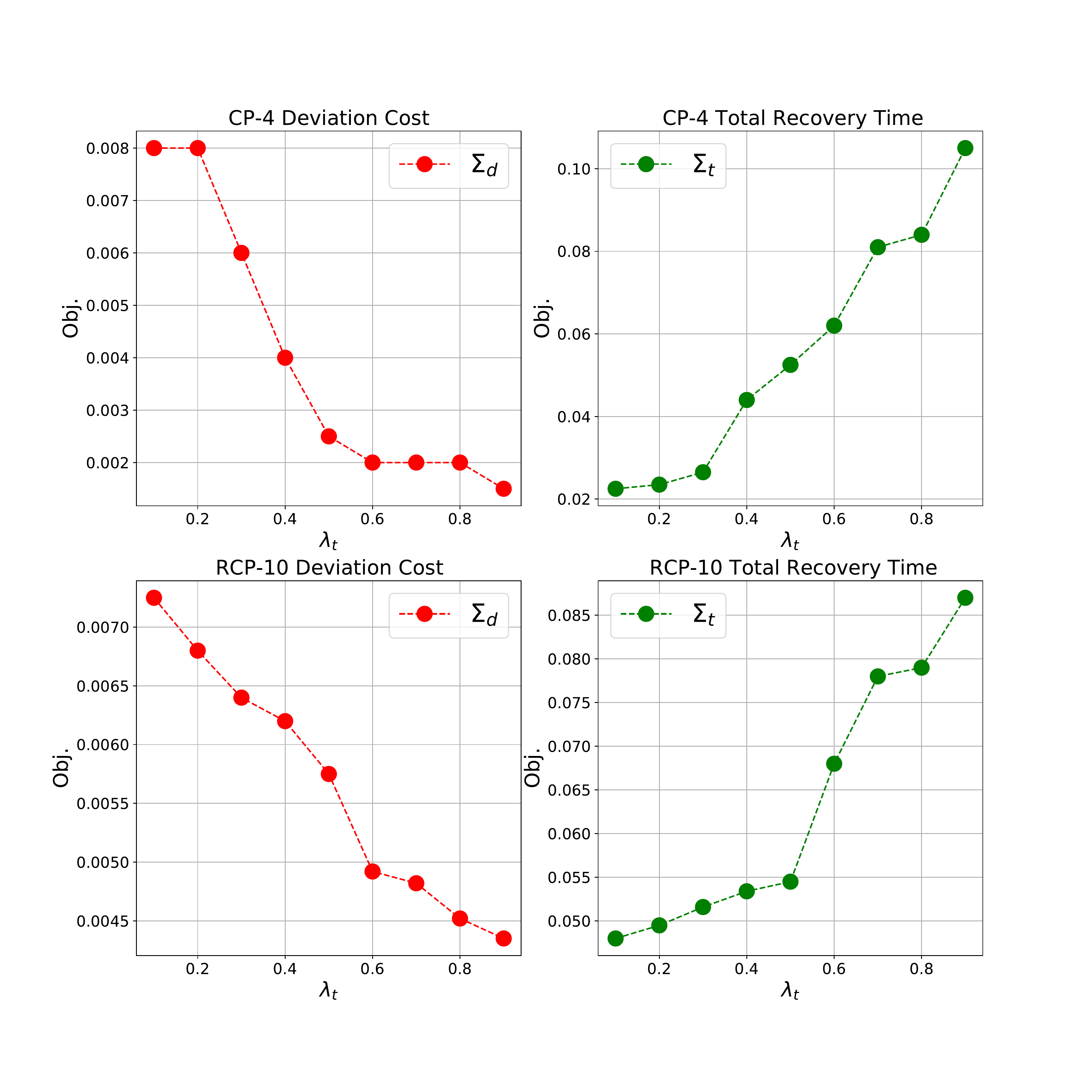}\label{scp10_6}}\\
	\caption{Sensitivity analysis on the preference weight $\lambda_t$ in the objective function Eq. \eqref{obj}. For all figures, $\Sigma_d$ represents the total speed and heading deviation defined as $\sum\limits_{i \in \mathcal{A}}(1-q_i)^2 + \theta_i^2$ (in red) and $\Sigma_t$ represents the total recovery time defined as $\sum\limits_{i \in \mathcal{A}} t_i^2$ (in green). All instances tested here are CP-4 and RCP-10.}
	\label{fig:sensitivitys}
\end{figure}

This behaviour is confirmed in our numerical experiments. Specifically, we solve the instances CP-4 and RCP-10-3 for $\lambda_t,\lambda_f = 0.1, \ldots, 0.9$ in steps of size $0.1$, i.e. for a total of 9 values per instance. Each parameter is changed separately. All instances are solved using Algorithm \ref{algo:2d_rec}. The change in the total deviation $\Sigma_d$ and the total recovery time $\Sigma_t$ are summarized in Figures \eqref{fig:sensitivityf} and \eqref{fig:sensitivitys}. In Figure \eqref{fig:sensitivityf}, it shows the behaviour of CP-4 and RCP-10-3 when $\lambda_f$ is changed: using the proposed objective function, the decision-maker can control which manoeuvre is prioritized by scaling up or down the preference weights  $\lambda_f$ accordingly. Higher values of $\lambda_f$ will minimize the total deviation, while smaller values of $\lambda_f$ will prioritize recovery time. We use $\lambda_f = 0.25$ in the numerical experiments presented in the remaining of the paper. We find that increasing $\lambda_f$ monotonically decreases the total deviation and monotonically increases the recovery time for all four instances tested. For $\lambda_t$ (in Figure \eqref{fig:sensitivitys}), higher values lead to the minimization of recovery time and bigger deviation while the increasing leads to smaller deviation in heading angle and speed while larger recovery time. We use $\lambda_t = 0.25$ in the numerical experiments presented in the remaining of the paper. For the parameter $w$, we are using the value $w=0.5$ based on the sensitivity analysis presented in \cite{dias2021disjunctive}.

\subsection{Performance of the Penalty-based Conflict Resolution and Trajectory Recovery Algorithm}
\label{perf_cont}

The performance of Algorithm \ref{algo:2d_rec} is reported in Table \ref{tab:CP-rec} for CP instances with 4 to 15 aircraft. The performance of the two-stage iterative algorithm is examined. The following tables represent the results for CP and RCP instances. In the header, there are four sections; the first is based on the instances containing $|\mathcal{A}|$ is the number of aircraft, and $n_c$ is the number of conflicts. The second group is based on the avoidance stage using the manoeuvre control variable: \textit{Obj.} is the objective function, $\deltatc$ is the optimality gap, Time(s) is the runtime in seconds, followed by total deviation in terms of $|1 -q_i|$, $|\theta_i|$ and $f_i$, respectively. The third group contains the results from the recovery stage using Model \ref{mod:onestepdetailed}. It is reported the average recovery time $\frac{1}{|A|}\sum\limits_{i \in \mathcal{A}}{t_i}$ among all aircraft $\min\limits_{i \in \mathcal{A}} t_i$. The fourth group is the overall performance of the algorithm corresponding to \textit{Iter.} with the number of iterations and overall runtime in seconds (s). The gap tolerance was 5\%. Table \ref{tab:CP-rec} summarises the results for CP instances and Table \ref{tab:RCP-rec} for the RCP instances with 10,20 and 30 aircraft. 

In the results, it is expected from those instances that all the aircraft are moving towards to centre. In terms of the objective function value the objective function is largely composed of $\sum\limits_{i \in \mathcal{A}}{f_i}$. For those instances, it is consistently observed that one aircraft can remain in its initial condition in each instance. This differs considerably from the results using naive avoidance, where the global solution would imply that all aircraft needs to perform some manoeuvres. In instances with ten or more aircraft, there are only slight deviations. Similarly, the variation in the heading angle is considered negligible for smaller instances. However, it does not increase proportionally with the number of conflicts, revealing that this manoeuvre is used sporadically.
Moreover, the current model results show that only a small amount of aircraft is required to be controlled, but the total deviation is more considerable. In terms of runtime, all instances with 12 aircraft can be solved within the time limit. In the recovery stage, the average recovery time does not increase with the number of aircraft (around 0.32 hours on average), with the minimum recovery time around 5 minutes while the maximum peaks at almost 1 hour. However, all the values of instances with 12 aircraft or more result in time out solutions. For recovery time, most instances can be solved within the time limit, but for the more prominent instances, the solution obtained in the avoidance stage time-out and the solution in the recovery stage. The algorithm only requires up to 3 iterations to obtain a solution within the convergence criterion used in terms of iterations. In comparison with the benchmark in Table \ref{tab:CP-sum}, we can conclude by the $TC$ value that the speed variation is quite smaller and the angle deviation is overall larger. This means that the recovery time is also considerably larger. This is mainly due to the trivial solutions where the main focus is to solve the conflict as it comes without any consideration relating to the recovery (as it is done in the iterative algorithm). Another observation is that the recovery strategies presented in the \Naive and in the \Greedy methods always provided a strategy where the recovery time was slightly higher for the latter and considerably higher for the former. This solidifies that using the cost of recovery as a preemptive measure contributes to better performance and usage of the airs space configuration.

\begin{sidewaystable}
\sidewaystablefn%
\begin{center}
\begin{minipage}{\textheight}
\caption{Summary of results for 2D CP instances with a speed control range of $[-6\%,+3\%]$ and a heading control range of $[-30^\circ,+30^\circ]$. All runtimes (Time) are reported in seconds while $\frac{1}{|A|}\sum\limits_{i \in \mathcal{A}}{t_i}$ is in hours and Gap(\%) expressed in percentage using the formulation applied in this papers.}\label{tab:CP-rec}
\resizebox{1.0\columnwidth}{!}{%
\begin{tabular}{ll llllll lll llll}
\toprule
&&\multicolumn{6}{l}{Avoidance} & \multicolumn{3}{l}{Recovery} & \multicolumn{4}{l}{Penalty-Based Algorithm}\\
\cmidrule(l){3-8} \cmidrule(l){9-11} \cmidrule(l){12-15}
$|A|$ & $|n_c|$ & Obj. & Gap (\%) & Time (s) & $\sum\limits_{i \in \mathcal{A}}{|1-q_i|}$  & $\sum\limits_{i \in \mathcal{A}}{|\theta_i}|$ & $\sum\limits_{i \in \mathcal{A}}{f_i}$ & Obj. & Time (s) & $\frac{1}{|A|}\sum\limits_{i \in \mathcal{A}}{t_i}$ & $\tc^0$ & $\tc^n$ & Iter. & Time(s)\\
\midrule
CP-4	&	6	&	3.00	&	1.79E-5	&	0.02	&	0.001	&	0.07	&	3 	&   3.39E-4	&	0.23 &	0.36	& 3.07 & 2.83 &	2	&	0.25 \\
CP-5	&	10	&	4.00	&	9.93E-5	&	0.02	&	0.001	&	0.11	&	4   &   1.19E-4	&	0.30 &	0.11	& 4.11 & 3.91 &	2	&	0.31	\\
CP-6	&	15	&	5.00	&	9.94E-5	&	0.03	&	0.002	&	0.08	&	5	&   6.37E-4	&	0.66 &	0.35	& 5.08 & 4.57 & 3	&	0.69	\\
CP-7	&	21	&	6.00	&	9.88E-5	&	0.11	&	0.001	&	0.12	&	6 	&   3.45E-4	&	3.85 &	0.26	& 6.20 & 5.58 &	3	&	3.96	\\
CP-8	&	28	&	7.00	&	9.94E-5	&	0.23	&	0.009	&	0.16	&	7 	&   9.71E-4	&	1.86 &	0.30	& 7.17 & 6.74 &	3	&	2.09	\\
CP-9	&	36	&	8.00	&	1.00E-4	&	2.89	&	0.000	&	0.22	&	8 	&   6.38E-4	&	15.1 &	0.20	& 8.22 & 7.40 &	3	&	18.0	\\
CP-10	&	45	&	9.01	&	9.99E-5	&	15.4	&	0.058	&	0.22	&	9 	&   1.63E-3	&	10.6 &	0.30	& 9.29 & 8.36 &	3	&	26.1	\\
CP-11	&	55	&	10.01	&	1.00E-4	&	59.2	&	0.010	&	0.30	&	10	&	1.56E-3	&	29.6 &	0.27	& 10.3 & 9.60 &	2	&	355 	\\
CP-12	&	66	&	11.01	&	1.00E-4	&	230	    &	0.000	&	0.36	&	11	&	2.29E-3	&	34.2 &	0.27	& 11.3 & 10.4 &	2	&	264 	\\
CP-13	&	78	&	12.01	&	3.59E-4	&	300	    &	0.050	&	0.36	&	12	&	2.05E-3	&	55.9 &	0.24	& 12.4 & 11.4 &	2	&	355 	\\
CP-14	&	91	&	13.01	&	5.83E-4	&	300	    &	0.010	&	0.40	&	13	&	3.27E-3	&	52.3 &	0.27	& 13.4 & 12.3 &	2	&	352 	\\
CP-15	&	105	&	14.02	&	6.80E-4	&	300	    &	0.080	&	0.46	&	14	&	3.22E-3	&	300	 &	0.25	& 14.5 & 13.3 &	2	&	600 	\\
\bottomrule
\end{tabular}
}
\end{minipage}
\end{center}
\end{sidewaystable}

\begin{sidewaystable}
\sidewaystablefn%
\begin{center}
\begin{minipage}{\textheight}
\caption{Summary of results for 2D CP instances with a speed control range of $[-6\%,+3\%]$ and a heading control range of $[-30^\circ,+30^\circ]$. All runtimes (Time) are reported in seconds while $\frac{1}{|A|}\sum\limits_{i \in \mathcal{A}}{t_i}$ is in hours and Gap(\%) expressed in percentage using the formulation applied in this papers.}\label{tab:RCP-rec}
\resizebox{1.0\columnwidth}{!}{%
\begin{tabular}{ll llllll lll llll}
\toprule
&&\multicolumn{6}{l}{Avoidance} & \multicolumn{3}{l}{Recovery} & \multicolumn{4}{l}{Penalty-Based Algorithm}\\
\cmidrule(l){3-8} \cmidrule(l){9-11} \cmidrule(l){12-15}
$|A|$ & $|n_c|$ & Obj. & Gap (\%) & Time (s) & $\sum\limits_{i \in \mathcal{A}}{|1-q_i|}$  & $\sum\limits_{i \in \mathcal{A}}{|\theta_i}|$ & $\sum\limits_{i \in \mathcal{A}}{f_i}$ & Obj. & Time (s) & $\frac{1}{|A|}\sum\limits_{i \in \mathcal{A}}{t_i}$ & $\tc^0$ & $\tc^n$ & Iter. & Time(s)\\
\midrule
RCP-10 & 3.10 (1.5) & 2.16 (0.93) & 3.4E-5 (3.5E-5) & 0.01 (0.00) & 0.02 (0.03) & 0.02 (0.02) & 2.16 (0.93) & 5.28E-5 (1.7E-4) & 2.14 (3.04) & 2.42 (2.26) & 2.6 (1.60) & 5.04 (1.01) & 4.73 (0.95) & 2.16 (3.05) \\
RCP-20 & 13.4 (3.4) & 6.66 (1.19) & 8.2E-5 (2.4E-5) & 0.33 (0.31) & 0.10 (0.08) & 0.12 (0.06) &  6.67 (1.19) & 4.86E-4 (4.1E-4) & 259  (68.8) & 6.96 (1.87) & 3.4 (0.50) & 45.4 (1.43) & 42.6 (1.34) & 259  (61.1)\\
RCP-30 & 33.7 (5.7) & 12.6 (1.46) & 9.1E-5 (1.6E-5) & 6.46 (9.02) & 0.39 (0.18) & 0.40 (0.15) & 12.6 (1.47) & 8.32E-3 (1.0E-2) & 300         & 10.6 (2.53) & 1.1 (0.21) & 159  (2.21) & 150  (2.08) & 600  \\ 
\bottomrule
\end{tabular}
}
\end{minipage}
\end{center}
\end{sidewaystable}

\begin{sidewaystable}
\sidewaystablefn%
\begin{center}
\begin{minipage}{\textheight}
\caption{Summary of Comparison between the Penalty-Based Algorithm and the Exact-Naive and Greedy-Naive presented in \cite{dias2020two} using the $TC$ as comparison key for the CP Instances.}\label{tab:CP-sum}
\resizebox{1.0\columnwidth}{!}{%
\begin{tabular}{l lll lll lll}
\toprule
&\multicolumn{3}{l}{Penalty-Based} & \multicolumn{3}{l}{\Naive} & \multicolumn{3}{l}{\Greedy}\\
\cmidrule(l){2-4} \cmidrule(l){5-7} \cmidrule(l){8-10}
$|A|$ & $\sum\limits_{i \in A} TC_i$ & $\frac{1}{|A|}\sum\limits_{i \in A} TC_i$ & Time & $\sum\limits_{i \in A} TC_i$ & $\frac{1}{|A|}\sum\limits_{i \in A} TC_i$ & Time & $\sum\limits_{i \in A} TC_i$ & $\frac{1}{|A|}\sum\limits_{i \in A} TC_i$ & Time \\
\midrule
4	&	2.83 	&	0.71	&	0.25	&	3.07	&	0.77	&	0.21	&	3.12	&	0.78	&	0.11	\\
5	&	3.91	&	0.78	&	0.31	&	4.11	&	0.82	&	0.26	&	4.56	&	0.91	&	0.13	\\
6	&	4.57	&	0.76	&	0.69	&	5.08	&	0.85	&	0.59	&	5.43	&	0.92	&	0.25	\\
7	&	5.58	&	0.80	&	3.96	&	6.20	&	0.89	&	3.37	&	6.91	&	0.98	&	1.51	\\
8	&	6.74	&	0.84	&	2.09	&	7.17	&	0.90	&	1.78	&	7.42	&	0.92	&	0.82	\\
9	&	7.40	&	0.82	&	18.0	&	8.22	&	0.91	&	15.3	&	8.71	&	0.96	&	6.85	\\
10	&	8.36	&	0.84	&	26.1	&	9.29	&	0.93	&	22.1	&	9.35	&	0.93	&	9.92	\\
11	&	9.60	&	0.87	&	355	    &	10.3	&	0.94	&	301 	&	10.8	&	0.98	&	56.1	    \\
12	&	10.40	&	0.87	&	264	    &	11.3	&	0.94	&	224	    &	12.1	&	1.08	&	97.1	    \\
13	&	11.40	&	0.88	&	355	    &	12.4	&	0.95	&	301	    &	13.9	&	1.06	&	135	    \\
14	&	12.30	&	0.88	&	352	    &	13.4	&	0.96	&	299	    &	14.5	&	1.03	&	164	    \\
15	&	13.30	&	0.89	&	600	    &  	14.5	&	0.97	&	600	    &	15.3	&	1.02	&	270 	\\
\bottomrule
\end{tabular}
}
\end{minipage}
\end{center}
\end{sidewaystable}

\begin{sidewaystable}
\sidewaystablefn%
\begin{center}
\begin{minipage}{\textheight}
\caption{Summary of Comparison between the Penalty-Based Algorithm and the Exact-Naive and Greedy-Naive presented in \cite{dias2020two} using the $TC$ as comparison key for the RCP Instances}\label{tab:RCP-sum}
\resizebox{1.0\columnwidth}{!}{%
\begin{tabular}{l lll lll lll}
\toprule
&\multicolumn{3}{l}{Penalty-Based} & \multicolumn{3}{l}{Exact Naive} & \multicolumn{3}{l}{Greedy Naive}\\
\cmidrule(l){2-4} \cmidrule(l){5-7} \cmidrule(l){8-10}
$|A|$ & $\sum_{i \in A} TC_i$ & $\frac{1}{|A|}\sum\limits_{i \in A} TC_i$ & Time & $\sum\limits_{i \in A} TC_i$ & $\frac{1}{|A|}\sum\limits_{i \in A} TC_i$ & Time & $\sum_{i \in A} TC_i$ & $\frac{1}{|A|}\sum\limits_{i \in A} TC_i$ & Time \\
\midrule
RCP-10 & 4.73 (0.95) & 0.47 (0.08) & 2.16 (3.05) & 5.07 (1.01) & 0.50 (0.00) & 1.16 (2.25) & 7.41 (1.41) & 0.74 (0.01) & 0.78 (1.35)\\
RCP-20 & 42.6 (1.34) & 2.17 (0.06) & 250 (61.1)  & 45.4 (1.43) & 2.13 (0.36) & 127 (45.6)  & 63.5 (3.51) & 3.17 (0.02) & 2.35 (2.14) \\
RCP-30 & 150 (2.08)  & 5.0 (0.03)  & 600         & 159  (2.21) & 5.12 (0.83) & 356 (132)   & 204 (12.7) & 6.81 (0.31) & 5.71 (4.31) \\
\bottomrule
\end{tabular}
}
\end{minipage}
\end{center}
\end{sidewaystable}

In Table \ref{tab:RCP-rec}, the results for RCP instances are presented. For the avoidance stage, the objective function reflects the number of aircraft required to be separated. As indicated by $f_i$, 2.16 for RCP-10, 6.66 for RCP-20 and 12.6 for RCP-30, which reflect the value obtained by $\sum\limits_{i \in \mathcal{A}}{f_i}$. Comparing the two values, it is noticeable that a large component of the objective function is solely based on that variable. In terms of speed, most aircraft do not perform any relatively large deviation. The values are close to the nominal value, with only 0.02 for RCP-10, 0.10 for RCP-20 and 0.39 for RCP-30. For heading changes, the values are small, reflecting that most aircraft do not perform any deviation in heading either: 0.02 for RCP-10, 0.12 for RCP-20 and 0.40 for RCP-30. However, compared to the heading deviation obtained using the complex number formulation, it is clear that heading deviation is considerably smaller and that most of the total deviation is caused by speed changes. The opposite behaviour is observed here. The runtime for those instances is reasonably short, and instances with up to 30 aircraft can be solved in less than 10 s. Compared with the complex number formulation results, the \Dis presented gives solutions balanced with equal contribution from speed and heading control. These current results state that speed and heading control are kept to their minimum while most of the control is set around whether certain aircraft need to be altered. Finally, the optimality gap is negligible in all instances. In the recovery stage, the runtime increases considerably with the number of aircraft given that the number of alternative routes to solve increases drastically: up to 2.14 s for RCP-10, 259 s for RCP-20 and time out for 100\% of the instances under 5 minutes of the time limit. In terms of iteration, RCP-10 instances and RCP-20 instances require up to 4 iterations to achieve convergence while solving RCP-30 instances; only up to two iterations are executed due to the time limit. In comparison with the benchmark methods, the table \ref{tab:RCP-sum}, the results showed that via analyzing the $TC$ values for the benchmarks, the recovery time in the naive recovery is significantly larger than the values obtained in the algorithm of this paper and in the greedy recovery is slightly bigger but only by a small margin. On the other hand, the runtime observed by the greedy recovery is relatively small, which highlights the computational trade-off at stake.

\section{Conclusion}
\label{con-rec}

The findings are summarised in Section \ref{find_chap6} and future research directions are discussed in Section \ref{future_chap6}.

\subsection{Summary of Findings}
\label{find_chap6}

A new mixed-integer formulation and a penalty-based algorithm for ACRP with trajectory recovery are proposed. The performance of this approach revealed that by echoing the cost associated with trajectory recovery, the avoidance could be manipulated to minimize pre-emptively the overall cost of ACRP. On average, in a couple of iterations, it showed that most instances could force the aircraft to have more significant deviations in the avoidance stage but with an earlier recovery time in the recovery stage. In this case, a trade-off between avoidance and recovery is observed. Another advantage of this addition is the concept of stability of the solutions throughout iterations. One of the main issues that happen by having an iterative algorithm is the variation in the profile of the solution. In each iteration,  the algorithm provides a different group of aircraft as a possible candidate to recover, and each solution has a different deviation cost, which can be an improvement. As expected, the deviation cost is higher with a higher deviation angle, but the number of aircraft manoeuvring is small. In this situation, adding the $f_i$ variable introduces an extra hurdle in modifying the status of an aircraft. With multiple iterations, it can be expected that less deviation in the whole set of aircraft will be observed, although a more significant deviation in individual aircraft. This creates stability in the solution. In the point of view of air traffic controllers, solutions where there are many aircraft manoeuvring lead to higher workload and are not desirable or implementable in real applications. 

The comparison between the naive and the iterative approaches shows that the latter can outperform the former in a set of criteria. First, the naive approach shows that the total cost is relatively more significant because there is no compensation for the recovery costs while the avoidance is processed. In this way, the naive approach overcompensates the avoidance manoeuvre and ignores less trivial solutions where small avoidance could have been taken. On the other hand, using the iterative approach starts with a solution similar to the naive approach and iteratively reaches towards more balanced and non-trivial solutions where we observed that a smaller avoidance angle is obtained; hence a later recovery time also be applied. Second, the naive approach does not incorporate the concept of manoeuvre control. There is a lack of control of the minimum number of aircraft trajectories in that scenario. The naive approach follows the trivial optimization problem, whereas a flaw in the original formulation favours a symmetrical solution where all aircraft are required to perform some manoeuvre with a near similar deviation. Therefore, all aircraft in such instances will alter their trajectories. In the iterative approach, this assumption is optimized where only a subset of aircraft change their speed and angle. This solution is a bit more realistic, given that in actual airspace occasions, fewer aircraft will be impacted by eventual conflict. Another aspect of this decision choice is that the aircraft required to alter their trajectory has significant manoeuvres instead of the solution using the naive approach where most solutions are pretty small and even hard to be detected. Finally, the iterative approach has more time leverages as continuous variables, allowing for a broader range of possible trajectories as the avoidance angle. The main setback of the iterative approach is that the convergence of such a model is imposed via the reduction of the total cost, which might lead to many iterations if the convergence rate is low. Also, the optimality condition is not guaranteed in both approaches, although the iterative is the straightforward best approach given all the reasons explained throughout this paper.

Alternatively, an adaptation of the algorithm in \cite{dias2021disjunctive} is proposed to solve the avoidance stage considering speed, heading control and manoeuvre control as decision variables. This is incorporated into an iterative two-stage algorithm that incorporates the projected cost of recovering the aircraft into the avoidance stage. A relatively small number of iterations showed that it could increase the overall deviation in the avoidance pre-emptively to reduce the recovery cost and, ultimately, the total cost. In the numerical experiments, the performance of the proposed algorithms shows that the number of aircraft required to be altered is reduced. Therefore, the deviation is more significant in heading angles for those required to be changed. However, the recovery time is shorter compared to the complex number formulation. In comparison with the \Dis results, it also shows that while the former provides solutions that are more balanced in terms of manoeuvres, they are also more invasive because it affects the whole set of aircraft. The algorithm present in this study presents an alternative solution where fewer aircraft are controlled will higher deviations, suggesting that this behaviour improves the total cost.

\subsection{Future Research and Perspectives}
\label{future_chap6}

The most significant limitation in the iterative two-stage algorithm is conflict avoidance and trajectory recovery decomposition. Although this is a common practice in mathematical programming, there is no guarantee of the quality of the solution. Ideally, a global solution should be created by jointly optimizing both stages in a unified formulation. This is heavily challenged by the non-linearity and complexity of such formulations. Therefore, further research is needed to model this problem more simply and efficiently. In addition, the cost of recovery is projected into the avoidance stage throughout the iterations. Alternatively, stochastic optimization methods can be used to determine the expected cost of recovery manoeuvres already incorporated in the avoidance stage. Finally, the discretization of any variable is a limitation and modelling recovery time as a continuous variable may help reduce the total cost of trajectories.

Further testing is required to fully assess the impact of the control variables on solution quality. Notably, quantifying the impact of first stage formulations with non-discretized heading changes and second stage formulation continuously is critical to evaluate the cost of manoeuvre discretization. Both formulations presented here are deterministic methods where no uncertainty is considered. Future research will focus on introducing uncertainty with the proposed formulation to incorporate trajectory prediction uncertainty capable of accounting for the expected cost of recovery at the conflict avoidance stage and develop robust conflict resolution algorithms to account for the uncertainty of aircraft trajectory prediction when planning.

%\begin{appendices}

%\section{Section title of first appendix}\label{secA1}

%An appendix contains supplementary information that is not an essential part of the text itself but which may be helpful in providing a more comprehensive understanding of the research problem or it is information that is too cumbersome to be included in the body of the paper.




%\end{appendices}


\bibliography{sn-bibliography}% common bib file


\end{document}